\documentclass{elsart}
\usepackage{color}
\usepackage{graphicx}
\usepackage{amsmath}
\usepackage{amssymb}
\newcommand{\nb}{n_{\bullet}}
\newcommand{\nw}{n_{\circ}}
\def\cT{\mathcal{T}}
\def\cA{\mathcal{A}}
\def\cB{\mathcal{B}}
\def\cC{\mathcal{C}}
\def\cO{\mathcal{O}}
\def\cE{\mathcal{E}}
\def\cP{\mathcal{P}}
\def\eleft{\epsilon_{\mathrm{left}}}
\def\eright{\epsilon_{\mathrm{right}}}
\def\prout{\mathcal{P}_r^{\mathrm{out}}}
\def\prin{\mathcal{P}_r^{\mathrm{in}}}
\def\pbout{\mathcal{P}_b^{\mathrm{out}}}
\def\pbin{\mathcal{P}_b^{\mathrm{in}}}

\newtheorem{theorem}{Theorem}
\newtheorem{proposition}{Proposition}%[theorem]
\newtheorem{corollary}{Corollary}%[theorem]
\newtheorem{lemma}{Lemma}%[theorem]

\newenvironment{proof}{\begin{pf}}{\qed\end{pf}}

\begin{document}
\begin{frontmatter}

\title{Transversal structures on triangulations: a combinatorial study and straight-line drawings}
\author{\'Eric Fusy}
\address{Algorithms Project, INRIA Rocquencourt and LIX, \'Ecole Polytechnique}
\ead{eric.fusy@inria.fr}

\begin{abstract}
 This article focuses on a combinatorial structure specific to triangulated plane graphs 
 with quadrangular outer face and no separating triangle, which are called irreducible triangulations.
 The structure has been introduced by Xin He under the name of regular edge-labelling
 and consists of two bipolar orientations that are transversal. For this reason, the terminology  used here 
 is that of transversal structures. The main results obtained in the article are a bijection between
 irreducible triangulations and ternary trees, and a straight-line drawing algorithm 
 for irreducible triangulations. For a random irreducible triangulation with $n$ vertices, the 
 grid size of the drawing is asymptotically with high probability $11n/27\times 11n/27$
 up to an additive error of $\cO(\sqrt{n})$. In contrast, the 
best previously known algorithm for 
 these triangulations only guarantees a grid size 
$(\lceil n/2\rceil -1)\times \lfloor n/2\rfloor$. 
\end{abstract} 
\end{frontmatter}

\section{Introduction}
A \emph{plane graph}, or planar map, 
is a connected graph embedded in the plane without
edge-crossings, and considered up to orientation-preserving homeomorphism. 
Many drawing algorithms for plane
graphs~\cite{To,Fel,deFrOss,Fe01,He93,Kant,Sch} endow
the graph with a particular structure, from which coordinates are assigned
to vertices in a natural way. For
example, triangulations, i.e., plane graphs with only triangular faces, 
are characterized by the fact that their inner edges can
 be partitioned into three spanning trees with specific incidence relations,
 the so-called Schnyder woods~\cite{Sch}. These
spanning trees yield a natural method to assign coordinates to a vertex,
by counting the number of 
faces in specific regions. 
Placing the vertices accordingly on an integer
 grid and linking adjacent vertices by
segments yields a straight-line drawing algorithm, which can be
refined to produce a drawing on an integer grid 
$(n-2)\times(n-2)$, with $n$ the number of vertices, see~\cite{Fel,Sch2}.

In this article we focus on so-called \emph{irreducible triangulations}, 
which are plane graphs with
quadrangular outer face, triangular inner faces, and no separating 
triangle, i.e., each 3-cycle is the
boundary of a face. 
Irreducible triangulations form an important 
class of triangulations, as they are closely related to 4-connected 
triangulations. In addition, as discussed in~\cite{Tri4},  many plane graphs,
including bipartite plane graphs and 4-connected plane graphs with at least 4 outer vertices,
can be triangulated into
an irreducible triangulation. There exist more compact straight-line
drawing algorithms for irreducible triangulations~\cite{Xin,Miura},
the size of the grid being guaranteed to be $(\lceil n/2\rceil
-1)\times \lfloor n/2\rfloor$ in the worst case. 
By investigating a bijection with ternary trees, we have observed that 
each irreducible triangulation $T$  can be endowed with a so-called
\emph{transversal structure}, which can be summarized as follows. 
Calling $W$, $N$, $E$, $S$ (like West, North, East, 
South) the four outer vertices of $T$ in clockwise order, the inner edges 
of $T$ can be oriented and partitioned into two sets: red edges that ``flow'' 
from $S$ to $N$, and blue edges that ``flow'' from $W$ to $E$.  
%As we will see,
%the structure can also be seen as a pair of bipolar
%orientations (one in each color) that are transversal. 
As we learned after
completing a first draft of this paper, X. He~\cite{He93} has defined the same
structure under the name of regular edge-labelling, and derived 
 a nice algorithm of
 rectangular-dual 
drawing, which has been recently applied to the theory of 
cartograms~\cite{Spe,SpeVa}. 
We give two equivalent definitions of 
transversal structures in Section~\ref{sec:tra}: one without orientations 
called a transversal edge-partition, and one with orientations called  a 
transversal pair of bipolar orientations (which corresponds 
to the regular edge-labelling
of X. He). Transversal structures characterize 
irreducible triangulations  in the same way as Schnyder Woods 
characterize triangulations, and they share similar combinatorial properties. 
In particular, we show in Section~\ref{sec:lattice} that the set of 
transversal structures of an irreducible triangulation  
is a distributive lattice, and that the ``flip'' operation has a 
simple geometric interpretation (Theorem~\ref{theo:lattice}). 

The transversal structure at the bottom of the lattice, called \emph{minimal},
 has a strong combinatorial role, as it
allows us to establish a bijection between ternary trees and 
irreducible triangulations.
The bijection, called \emph{closure mapping}, 
is described in Section~\ref{sec:bij}. 
The mapping from ternary trees to irreducible triangulations 
relies on ``closure operations'', as introduced by G. 
Schaeffer in his PhD~\cite{S-these}, 
see also~\cite{Fu,PS03}. This bijection has in fact brought
about our discovery of transversal edge-partitions, as
 a  natural \emph{edge-bicoloration} of a ternary tree is mapped 
to the minimal transversal edge-partition of the associated triangulation 
(similary, the bijection of~\cite{PS03} maps the structure of 
Schnyder woods). Classical algorithmic applications of bijections 
between trees and plane graphs are random
generation ---as detailed in~\cite{S-these,Sc99}--- 
and encoding algorithms for plane graphs ---as detailed 
in~\cite{Fu,PS03}--- with application to mesh compression in computational geometry.
The closure mapping presented in this article yields linear time procedures for random
sampling (under a fixed-size uniform distribution) and optimal  encoding (in the information theoretic sense) of \emph{4-connected triangulations}.
These algorithms are
described in the thesis of the author~\cite{Fu-these}. 
The focus in this article (besides the graph drawing algorithm) is on the application to counting;  the bijection yields a
combinatorial way to
enumerate rooted 4-connected triangulations, which were already
counted by Tutte in~\cite{Tut} using algebraic methods.

In Section~\ref{sec:draw}, we derive from transversal structures a 
straight-line drawing algorithm for irreducible triangulations. 
In a similar way as algorithms using Schnyder Woods~\cite{Fel,Sch}, 
the drawing is obtained by using face-counting operations. 
Our algorithm outputs a straight line drawing on an integer grid 
of half-perimeter $n-1$ if the triangulation has $n$
vertices (Theorem~\ref{theo:draw}).  
This is to be compared with previous algorithms 
for irreducible triangulations
by He~\cite{Xin}
and Miura \emph{et al}~\cite{Miura}. The latter produces a grid
 $(\lceil n/2\rceil -1)\times \lfloor n/2\rfloor$; the 
half-perimeter is also $n-1$,
but the aspect ratio of the outer face is better. However,
the algorithms of~\cite{Xin}
and~\cite{Miura} rely on a particular order to treat the vertices,
called \emph{canonical ordering}, and a step of
coordinate-shifting makes them difficult to implement and to carry out
by hand. In contrast, our algorithm can readily be performed on
a piece of paper, because the coordinates of the vertices are computed
independently with simple face-counting operations.

Furthermore, some coordinate-deletions can be performed on 
the drawing obtained
using the face-counting algorithm, with the effect of 
reducing the size of the grid (Theorem~\ref{theo:compactdraw}). 
For an irreducible
triangulation with $n$ vertices taken uniformly at 
random and endowed with its minimal transversal
structure (for the distributive lattice), we show in 
Section~\ref{sec:proofreduc} (Theorem~\ref{prop:reduction}) that the
size of the grid after coordinate-deletions is asymptotically 
with high probability
$11n/27\times 11n/27$ up to an additive error of order 
$\sqrt{n}$. 
Compared to~\cite{Xin}
and~\cite{Miura}, we
do not improve on the size of the grid in the worst case, but we improve
asymptotically with high probability by a reduction-factor 
$27/22$ on the width and height of the
grid, see Figure~\ref{fig:big} 
for an example. The proof of the grid size $11n/27\times 11n/27$ 
makes use of several ingredients:
a combinatorial interpretation of coordinate-deletions, the bijection with ternary trees, and modern tools of analytic 
combinatorics such as the quasi-power theorem~\cite{fla}. 

The following diagram summarizes the connections between the different 
combinatorial structures and the role they play for the drawing algorithms. 

\begin{center}
\includegraphics[width=10cm]{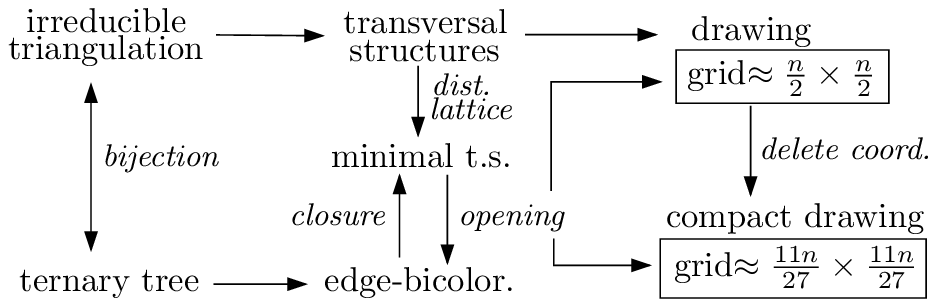}
\end{center}

\section{Transversal structures: definitions}
In this section we give two definitions of transversal structures, 
one without orientations called
transversal edge-partition and one with orientations called 
transversal pair of bipolar orientations.
As we will prove in Proposition~\ref{prop:lattice}, the two 
definitions are in fact equivalent, i.e., the 
additional information given by the orientations is redundant.
The definition without orientations 
is more convenient for the combinatorial study (lattice property 
and bijection with ternary trees), 
while the definition with orientations
fits better to describe  the straight-line 
drawing algorithm in Section~\ref{sec:draw}.

\begin{figure}[h]
\begin{center}
\includegraphics{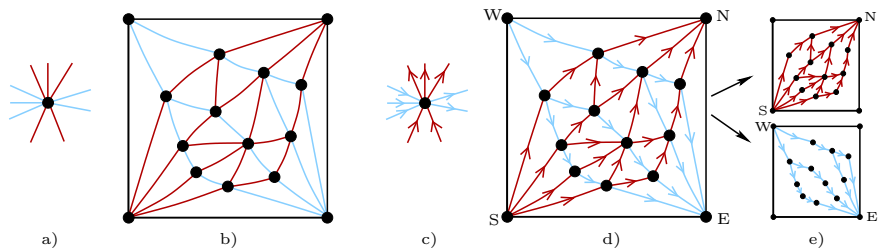}
\end{center}
\caption{Transversal edge-partition: local condition (a) 
and a complete 
example (b). In parallel, transversal pair of bipolar 
orientations: 
local condition (c) and a complete example (d) with the two induced
bipolar orientations (e).}
\label{fig:Example}
\end{figure}

\label{sec:tra}
\subsection{Transversal edge-partition}
\label{sec:tra1}
Let $T$ be an irreducible triangulation. Edges and vertices 
of $T$ are said to be \emph{inner} or \emph{outer} whether they are incident
 to the outer face or not. A \emph{transversal edge-partition} of $T$ 
is a combinatorial structure where each inner edge of $T$ is given
a color ---red or blue--- such 
that the following conditions are satisfied.

\begin{description}
\item[C1 (Inner vertices):]
In clockwise order around each inner vertex $v$, the edges incident to $v$ 
form: a non empty interval of red edges,  a non empty interval of blue edges, 
a non empty 
interval of red edges, and a non empty interval of blue edges, see 
Figure~\ref{fig:Example}(a).
\item[C2 (Outer vertices):] 
Writing $a_1$, $a_2$, $a_3$, $a_4$ for the outer vertices of $T$ in 
clockwise order, all inner edges incident to $a_1$ and to $a_3$ are of one 
color and all 
inner edges incident to $a_2$ and to $a_4$ are of the other color.
\end{description}

An example of transversal edge-partition is illustrated in 
Figure~\ref{fig:Example}(b), 
where red edges are
\textsc{dark grey} and blue edges are \textsc{light grey} 
(the same convention will
be used for all figures).

\subsection{Transversal pair of bipolar orientations}
\label{sec:bip}
An orientation of a graph $G$ is said to be \emph{acyclic} 
if it has no oriented circuit (a circuit is an oriented simple cycle 
of edges). Given an acyclic orientation of $G$, a vertex 
 having no ingoing edge is called a \emph{source}, and a vertex having no
  outgoing edge is called a \emph{sink}. A \emph{bipolar orientation} is 
an acyclic orientation with a unique source, denoted $s$, and a unique sink, denoted $t$.
Such an orientation is also characterized by the fact that, for each 
vertex $v\neq\{s,t\}$, there exists an oriented path 
from $s$ to $t$ passing by $v$, see~\cite{Oss} for a detailed discussion. 
An important property of a bipolar orientation on a \emph{plane graph} 
is that the edges 
incident to a vertex $v\neq \{s,t\}$ are partitioned into a non-empty interval
of ingoing edges and a non-empty interval of outgoing edges; and, dually, each
face $f$ of $M$ has two particular vertices $s_f$ and $t_f$ such that
the boundary of $f$ consists of two non-empty oriented paths both going
from $s_f$ to $t_f$, the path with $f$ on its right (left) is called the \emph{left lateral path} 
  (\emph{right lateral path}, resp.) of $f$.

Let $T$ be an irreducible triangulation. Call $N$, $E$, 
$S$ and $W$ 
the four outer vertices of $T$ in clockwise order around 
the outer face. 
A \emph{transversal pair of bipolar orientations} of $T$ 
is a combinatorial structure
where each inner edge of $T$ is given a direction and a 
color ---red or blue--- such that the
following 
conditions are satisfied (see Figure~\ref{fig:Example}(d) for an example):

\begin{description}
\item[C1' (Inner vertices):]
 In clockwise order around each inner vertex $v$ of $T$, the  
edges incident to $v$ form: a non empty interval of outgoing 
red edges, a non empty interval 
of outgoing 
blue edges, a non empty interval of ingoing red edges, and a non empty 
interval of ingoing 
blue edges, see Figure~\ref{fig:Example}(c).
\item[C2' (Outer vertices):]
 All inner edges incident to $N$, $E$, $S$ and $W$ are 
respectively ingoing red, ingoing blue, outgoing red, and outgoing blue.
\end{description}

This structure is also considered in~\cite{He93,Kant} under the name of
\emph{regular edge labelling}.

\begin{proposition}
\label{prop:acyc}
The orientation of the edges given by a 
transversal pair of bipolar 
orientation is acyclic. The sources are $W$ and $S$, and the 
sinks are $E$ and $N$.
\end{proposition} 
\begin{proof}
Let $T$ be an irreducible triangulation  endowed with a 
transversal pair of bipolar orientations. Assume there exists 
a circuit, and consider a minimal one $\cC$,
 i.e., the interior of $\cC$ is not included in the interior of any 
other circuit. It is easy to check from Condition C2' that 
$\mathcal{C}$ is not the boundary of a face of $T$. Thus the interior of 
$\mathcal{C}$ contains at least one edge $e$ incident to a vertex $v$ of 
$\mathcal{C}$. Assume that $e$ is going out of $v$. Starting from $e$, 
it is always possible, when reaching a vertex inside $\mathcal{C}$, to go 
out of that vertex toward one of its neighbours. Indeed, a vertex inside 
$\mathcal{C}$ is an inner vertex of $T$, hence has positive outdegree 
according to Condition C1'. Thus, there is an oriented path starting 
from $e$, that either loops $\mathcal{C}$ into a circuit in the interior of
$\cC$ ---impossible by minimality of $\cC$--- or reaches 
$\mathcal{C}$ again ---impossible as a chordal path for $\cC$
would yield two smaller circuits. 
Thus, the orientation is acyclic. 
Finally Condition C1' ensures that no inner vertex of $T$ 
can be a source or a sink,
and Condition C2' ensures that $W$ and $S$ are sources, and $E$ and $N$ are
sinks.\hfill\phantom{1}
\end{proof}

The following corollary of Proposition~\ref{prop:acyc}, 
also proved in~\cite{He93}, explains the name of
transversal pair of bipolar orientations, see also Figure~\ref{fig:Example}(e).

\begin{corollary}[\cite{He93}]
\label{prop:bip}
Let $T$ be an irreducible triangulation  endowed with a 
transversal pair of bipolar 
orientation. Then the (oriented) red edges induce a bipolar orientation 
of the plane 
graph obtained from $T$ by removing $W$, $E$, and all non red edges. 
Similarly, the 
blue edges induce a bipolar orientation of the plane graph 
obtained from $T$ by deleting $S$, $N$, 
and all non blue 
edges.
\end{corollary}

%\noindent\textbf{Remark:} As transversal pairs of bipolar orientations are in 
%bijection with transversal edge-partitions of a given irreducible 
%triangulation $T$, we also call \emph{minimal} the transversal 
%pair of bipolar orientations associated with the minimal transversal 
%edge-partition of $T$.

\section{Lattice property of transversal edge-partitions}
\label{sec:lattice}

We investigate on the set $\mathcal{E}(T)$ of transversal 
edge-partitions of a fixed irreducible triangulation $T$. 
Kant and He~\cite{Kant} have shown
that $\cE(T)$ is not empty and that an element of $\cE(T)$ 
can be computed in linear time. 
In this section, we prove that $\cE(T)$ is a 
distributive lattice. This property is to be compared with 
the lattice property of other 
similar combinatorial structures on plane graphs, such as 
bipolar orientations~\cite{Men} and
Schnyder woods~\cite{Bre,Fellat}.  Our proof takes advantage of the 
property that the set of orientations of a plane graph with a prescribed 
outdegree for each vertex is a distributive lattice. To apply this result to 
transversal edge-partitions, we establish a bijection between $\mathcal{E}(T)$ 
and some orientations with prescribed vertex outdegrees 
on an associated plane graph, called the \emph{angular graph} of $T$.

\subsection{Lattice structure of $\alpha$-orientations on a plane graph}
Let us first recall the definition of a distributive lattice. 
A \emph{lattice} is a partially ordered set $(E,\leq)$ such that,
for each pair
$(x,y)$ of elements of $E$, there exists a unique element $x\wedge y$ 
and a unique element $x\vee y$ satisfying the conditions:
\begin{itemize}
\item
$x\wedge y\leq x$, $x\wedge y\leq y$, and $\forall z\in E$, 
$z\leq x$ and $z\leq y$ implies $z\leq x\wedge y$,
\item
$x\vee y\geq x$, $x\vee y\geq y$, and $\forall z\in E$, 
$z\geq x$ and $z\geq y$ implies $z\geq x\vee y$.
\end{itemize}
In other words, each pair admits a unique common lower element dominating
all other common lower elements, and the same holds with common upper elements.
The lattice is said to be \emph{distributive} if the operators $\wedge$
and $\vee$ are distributive with respect to each other, i.e., $\forall
(x,y,z)\in E$, $x\wedge(y\vee z)=(x\wedge y)\vee(x\wedge z)$ 
and $x\vee(y\wedge z)=(x\vee y)\wedge (x\vee z)$. The nice 
feature of distributive lattices is that,
in most cases, moving from an element of the lattice to a 
covering lower (upper) element  has a
simple geometric interpretation, which we informally call 
a flip (flop, respectively).
As we recall next, in the case of 
orientations of a plane graph with prescribed vertex outdegrees, the
flip (flop) operation consists in reversing a clockwise circuit 
(counter-clockwise circuit, respectively).

Given a plane graph $G=(V,E)$, and a
function $\alpha :V\rightarrow \mathbf{N}$, an \emph{$\alpha$-orientation} of
$G$ is an orientation of $G$ such that each vertex $v$ of
$G$ has outdegree $\alpha(v)$. An oriented circuit $\mathcal{C}$ 
is called 
\emph{essential} if it has no chordal path (a \emph{chordal path} 
is an oriented path of edges in the interior of $\mathcal{C}$ with the two 
extremities on $\mathcal{C}$). 
%The following 
%theorem appeared first in 
%the thesis of Ossona de Mendez~\cite{Men} and was rediscovered by 
%Felsner in~\cite{Fellat}.

\begin{theorem}[Ossona de Mendez~\cite{Men}, Felsner~\cite{Fellat}]
\label{theo:alphaorientation}
Given $G=(V,E)$ a plane graph and $\alpha :V\rightarrow \mathbf{N}$ a
function, the set of $\alpha$-orientations of $G$ is either empty or
is a distributive lattice. The flip operation consists in
reversing the orientation of an essential clockwise circuit.
\end{theorem}

\subsection{Bijection with orientations of the angular graph}
\label{sec:corresTransAlpha}
\begin{figure}
\begin{center}
\includegraphics{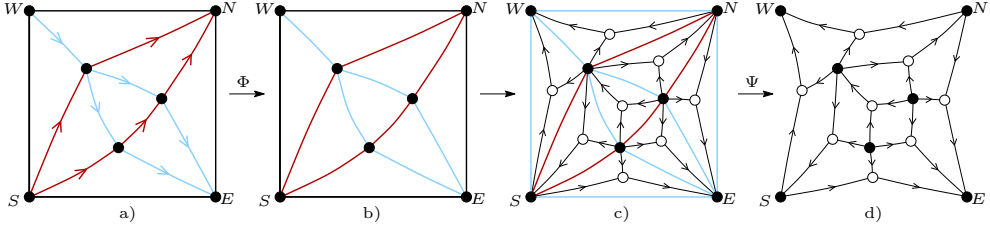}
\end{center}
\caption{Given an irreducible triangulation $T$  endowed
with a transversal pair of bipolar orientations $Z$ (Fig.a) 
and the induced transversal 
edge-partition $E=\Phi(Z)$ (Fig.b), construction of the 
angular graph $Q(T)$
and of the $\alpha_0$-orientation of $Q(T)$ image of $E$ by the mapping 
$\Psi$ (Fig.c and d).}
\label{fig:Angular}
\end{figure}

Let $T$ be an irreducible triangulation. The \emph{angular graph} of $T$ is the
bipartite plane graph $Q$ with vertex set $V_Q$  consisting of vertices 
and inner faces of $T$,
and edge set corresponding to the incidences between these vertices and
faces, see Figure~\ref{fig:Angular}. In the sequel, vertices of $Q$ 
corresponding
to vertices (inner faces) of $T$ are black (white, respectively), and the whole vertex set of $Q$
is denoted $V_Q$. 
The terminology is due
to the fact that each edge of $Q$ corresponds to an angle $(v,f)$ of $T$
(the angles incident to the outer face of $T$ are not considered).
Notice also that each inner face of $Q(T)$ is quadrangular. 
 
We consider the function $\alpha_0:V_Q\rightarrow \mathbf{N}$ specified 
as follows,
\begin{itemize}
\item
each black vertex $v$ of $Q(T)$  corresponding to an inner vertex of $T$ 
satisfies $\alpha_0(v)=4$,
\item
each white vertex $v$ of $Q(T)$ satisfies $\alpha_0(v)=1$,
\item
the outer black vertices satisfy $\alpha_0(N)=\alpha_0(S)=2,\ \alpha_0(W)=\alpha_0(E)=0$. 
\end{itemize}

\begin{proposition}
\label{prop:lattice}
Given $T$ an irreducible triangulation and $Q(T)$ the angular graph of $T$, 
the following sets
are in bijection:
\begin{itemize}
\item
transversal
pairs of bipolar orientations of $T$,
\item
transversal edge-partitions of $T$,
\item
$\alpha_0$-orientations of $Q(T)$.
\end{itemize}
%Hence
%Theorem~\ref{theo:alphaorientation} ensures that $\mathcal{B}$ and
%$\mathcal{E}$ have a lattice structure.
\end{proposition}
In particular, Proposition~\ref{prop:lattice} ensures that the 
definitions of transversal edge-partition and
of transversal pair of bipolar orientations are equivalent, i.e., 
the orientation of edges is 
a information once the colors are fixed.

Given $T$ an irreducible triangulation and $Q(T)$ the angular graph of $T$, 
we denote by $\cB$ the set
of transversal pairs of bipolar orientations of $T$, $\cE$ the set of 
transversal edge-partitions of $T$, 
and $\cO$ the set of $\alpha_0$-orientations of $Q(T)$.
To prove Proposition~\ref{prop:lattice}, we introduce two mappings $\Phi$
and $\Psi$ respectively from $\mathcal{B}$ to $\mathcal{E}$ and from
$\mathcal{E}$ to $\mathcal{O}$.

Given $Z\in \mathcal{B}$, $\Phi(Z)$ is simply the edge-bicoloration 
induced  by $Z$,
as illustrated in Figure~\ref{fig:Angular}(a)-(b). 
It is straightforward  that 
$\Phi(Z)\in \mathcal{E}$. 
In addition, $\Phi$ is
clearly injective. Indeed, starting from a transversal edge-partition, 
the directions
of edges of $Q(T)$ are assigned greedily  
so as to satisfy the local rules of a transversal pair of bipolar 
orientations. The fact that the propagation of edge directions is done 
without conflict is not straightforward; it will follow from the 
surjectivity of the mapping $\Phi$, to be proved later.

Given $X\in \mathcal{E}$, we define $\Psi(X)$ as the following
orientation of $Q(T)$. First, color blue the four outer
edges of $T$. Then, for each angle $(v,f)$ of $T$, orient the corresponding
edge of $Q(T)$ out of $v$ if $(v,f)$ is bicolored, 
and toward $v$ if $(v,f)$ is unicolored. 
(An angle $(v,f)$ of $T$ is called \emph{bicolored} 
if it is delimited by two edges of $T$ having different colors, 
and is called unicolored
otherwise.) 
Condition C1 implies that all inner black vertices of
$Q(T)$ have outdegree 4. In addition, Condition C2 and the fact 
that the four
outer edges of $T$ have been colored blue imply that $E$ and $W$ have
outdegree 0 and that $N$ and $S$ have outdegree 2.

The following lemma ensures that all white vertices have outdegree 1 
in $\Psi(X)$, so that 
$\Psi(X)$ is an $\alpha_0$-orientation.

\begin{lemma}
\label{lemma:bic}
Let $T$ be a plane graph with quadrangular outer face, triangular 
inner faces, and 
 endowed with a transversal edge-partition, the four outer
edges being additionally colored blue. 
Then there is no mono-colored inner face, i.e.,  each inner face of $T$ has two sides of one
color and one side of the other color. 
\end{lemma}
\begin{proof}
Let $\Lambda$ be the number of bicolored angles of $T$ and let $n$ be the 
number of inner vertices of
$T$. Condition C1 implies that there are $4n$
bicolored angle incident to an inner vertex of $T$.  
Condition C2 and
the fact that all outer edges are colored blue imply that two angles
incident to $N$ and two angles incident to $S$ are
bicolored. Hence, $\Lambda=4n+4$. 

Moreover, as $T$ has a 
quadrangular outer face and triangular inner faces, Euler's 
relation ensures that $T$ has 
$2n+2$ inner faces. For each
inner face, two cases can arise: either the three sides have the same
color, or two sides are of one color and one side is of the other
color. In the first (second) case, the face has 0 (2, respectively)
bicolored angles. As there are $2n+2$ inner faces and $\Lambda=4n+4$, 
the pigeonhole
principle implies that all inner faces have a contribution of 2 to the
number of bicolored angles, which concludes the proof.
\hfill\phantom{1}\end{proof}

Lemma~\ref{lemma:bic} ensures that $\Psi$ is a mapping from $\mathcal{E}$ to 
$\mathcal{O}$. In addition, it
is clear that $\Psi$ is injective. To prove Proposition~\ref{prop:lattice}, 
it remains to prove that $\Phi$
and $\Psi$ are surjective. As $\Phi$ is injective, it is sufficient 
to show that
$\Psi \circ \Phi$ is surjective. Thus, given $O\in \mathcal{O}$, we have to 
find $Z\in \mathcal{B}$ such that $\Psi \circ \Phi (Z)=O$. 
%To do that, we
%could simply greedily traverse the triangulation and propagate a
%coloration and orientation of its edges that is locally (around each
%vertex) compatible with the given $\alpha$-orientation $O$. However
%we have no guarantee that the same edge can not receive two different
%colors or orientations (thus in conflict) 
%during such a greedy traversal. It turns out
%that no such conflict can arise, as we prove next by
%traversing the vertices of the triangulation in a certain order.

\vspace{0.5cm}

\noindent\textbf{Computing the preimage of an 
$\alpha_0$-orientation.}
We now describe a method to compute a transversal pair of bipolar orientations
 $Z$ consistent with a given $\alpha_0$-orientation $O$, i.e., such that
  $\Psi\circ\Phi(Z)=O$. The algorithm makes use of a sweeping process
to orient and color the inner edges of $T$; a simple (i.e., not
self-intersecting) path $\cP$ of inner 
edges of $T$ going from $W$ to $E$ is maintained, the path moving progressively
toward the vertex $S$ (at the end, the path is $\cP=W\to S\to E$). We require
that the following invariants are satisfied throughout the sweeping process.
\begin{enumerate}
\item
Each vertex  of $\cP\backslash \{W,E\}$ has two outgoing edges on each side
of $\cP$ for the $\alpha_0$-orientation $O$.
\item
The inner edges of $T$ already oriented and colored are those on the left 
of $\cP$.
\item
Condition C1' holds around each inner vertex of $T$ on the left of $\cP$.
\item
A partial version of C1' 
holds around each vertex $v$ of $\cP\backslash\{W,E\}$.
The edges incident to $v$ on the left of $\cP$ form in clockwise order:
a possibly empty interval of ingoing blue edges, a non-empty interval of
outgoing red edges, and a possibly empty interval of outgoing blue edges.
\item
All edges already oriented and colored and incident to $N$, $E$, $S$, $W$ are
ingoing red, ingoing blue, outgoing red, and outgoing blue, respectively.
\item
The edges of $T$ already colored (and oriented) are consistent with the
$\alpha_0$-orientation $O$, i.e., for each angle $(v,f)$ delimited by two
edges of $T$ already colored, the corresponding edge of $Q(T)$ is going 
out of $v$ iff the angle is bicolored.  
\end{enumerate}

At first we need a technical result ensuring that Invariant~(1) is sufficient
for a path to be simple.

\begin{lemma}
\label{lem:simple}
Let $T$ be an irreducible triangulation and $Q(T)$ be the angular graph of $T$,
 endowed with an $\alpha_0$-orientation $O$. Let $\cP$ be a path of 
inner edges of $T$ from $W$ to $E$ such that each vertex of $\cP\backslash\{W,E\}$ has outdegree 2 on each side of $\cP$ for the $\alpha_0$-orientation.
Then the path $\cP$ is simple.
\end{lemma}
\begin{proof}
Assume that the path 
$\cP$ loops into a circuit; and consider an inclusion-minimal such
circuit $\cC=(v_0,v_1,\ldots,v_k=v_0)$. We define 
$\nb$, $\nw$ and $e$ as the numbers of black vertices (i.e., vertices of $T$), 
white vertices and edges of $Q(T)$ inside $\mathcal{C}$. 
As $T$ is triangulated and $\cC$ has length $k$, 
Euler's relation ensures that $T$ has  $2\nb+k-2$ faces inside $\cC$, i.e.,
$\nw=2\nb+k-2$.
Counting the edges of $Q(T)$ inside
 $\mathcal{C}$ according to their incident white vertex gives 
$(i):\ e\!=\!3\nw\!=\!6\nb\!+\!3k\!-\!6$. 
The edges of $Q(T)$
inside $\mathcal{C}$ can also be counted
according to their origin for the $\alpha_0$-orientation. 
Each vertex of $\cC$ ---except possibly the self-intersection vertex $v_0$---
  has outdegree 2 in the
interior of $\cC$ for the $\alpha_0$-orientation. Hence,
$(ii):\ e=4\nb+\nw+2k-2+\delta=6\nb+3k-4+\delta$, 
where $\delta\geq 0$ is the outdegree of $v_0$ inside $\cC$. 
Taking $(ii)-(i)$ yields $\delta=-2$, a contradiction. 
\hfill\phantom{1}  
\end{proof}

\begin{figure}
\begin{center}
\includegraphics{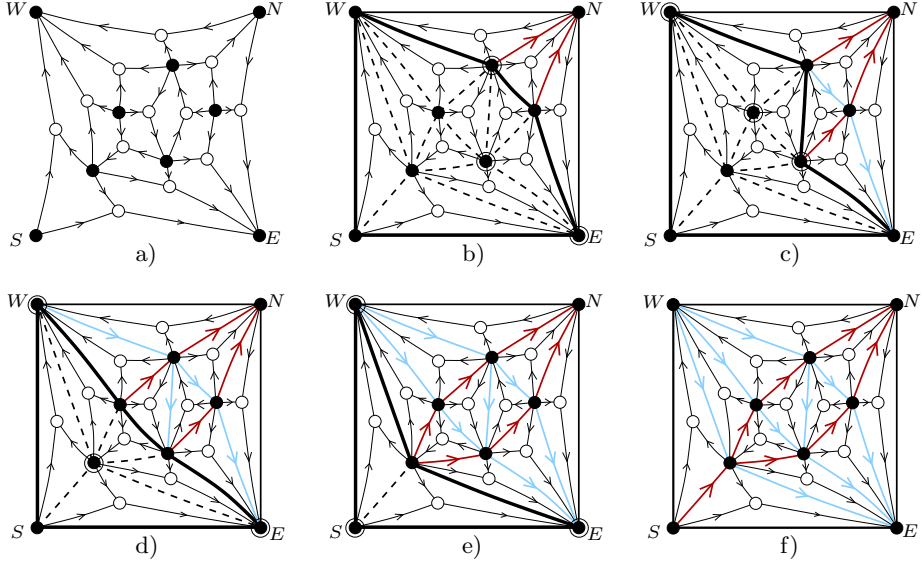}
\end{center}
\caption{The complete execution of the algorithm calculating the
preimage of an $\alpha_0$-orientation $O$. At each
step, the vertices of the matching path $P(v,v')$ are surrounded.}
\label{fig:AlgoBij}
\end{figure}

The path $\mathcal{P}$ is initialized with all neighbours
of $N$, from $W$ to $E$. 
In addition, all inner edges incident to
$N$ are initially colored red and directed toward $N$, see 
Figure~\ref{fig:AlgoBij}(b). The invariants (1)-to-(6) are clearly true
at the initial step.

\begin{figure}
\begin{center}
\includegraphics[width=10cm]{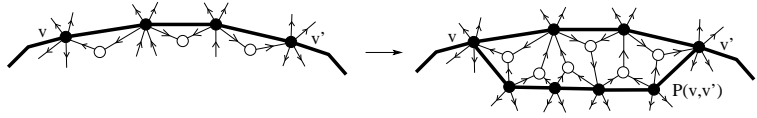}
\end{center}
\caption{An admissible pair $v,v'$ of vertices, and the matching path $P(v,v')$.}
\label{fig:upda_Path}
\end{figure}

\begin{figure}
\begin{center}
\includegraphics[width=10cm]{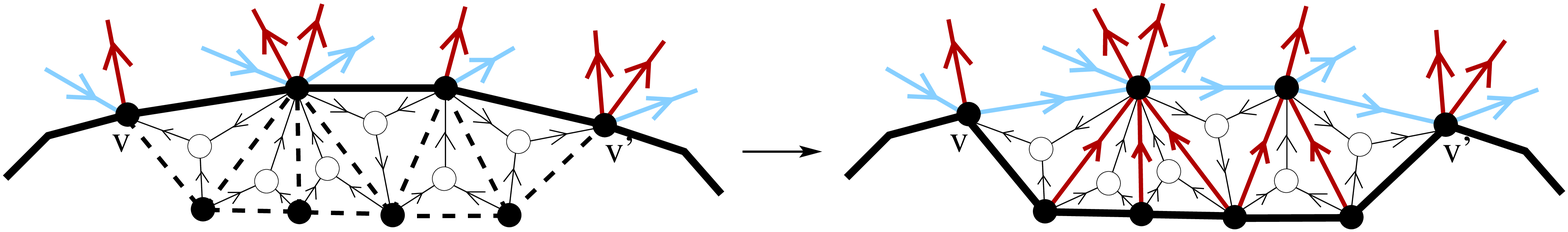}
\end{center}
\caption{The update step of the iterative algorithm to find the
preimage of an $\alpha_0$-orientation.}
\label{fig:upda}
\end{figure}

Let us introduce some terminology in order to describe the sweeping process. 
Thoughout the process, 
the vertices of $\cP$ are ordered from left to
right, with $W$ as leftmost and $E$ as rightmost vertex. Given $v,v'$ 
a pair of vertices on $\mathcal{P}$ ---with $v$ on the left of $v'$---
the part of $\cP$ going from $v$ to $v'$ is denoted by $[v,v']$.
For each vertex $w$ on $\cP$, let $f_1,\ldots,f_k$ be the sequence of 
faces of $T$ incident to $w$ on the right of $\cP$, taken in counterclockwise
order. The edge of $Q(T)$ associated to the angle $(w,f_1)$ (angle $(w,f_k)$)
is denoted by $\eleft(w)$ (by $\eright(w)$, respectively). A pair of vertices
$v,v'$ on $\cP$ ---with $v$ on the left of $v'$--- is called \emph{admissible}
if $\eright(v)$ is ingoing at $v$, $\eleft(v')$ is ingoing at $v'$, and for
each vertex $w\in [v,v']\backslash\{v,v'\}$, the edges $\eleft(w)$ and $\eright(w)$ are going out of $w$. Clearly an admissible pair always exists: take
a pair $v,v'$ of vertices on $\cP$ that satisfy \{$\eright(v)$ ingoing at $v$,
$\eleft(v')$ ingoing at $v'$\} and are closest possible. Notice that two  vertices $v,v'$ forming an admissible pair are not neighbours on $\cP$ (otherwise
the white vertex associated to the face on the right of $[v,v']$ would 
have outdegree $>1$). Let $w_0=v,w_1,\ldots,w_k,w_{k+1}=v'$ ($k\geq 1$) be the sequence
of vertices of $[v,v']$. The \emph{matching path} 
of $v,v'$ is the path of edges
of $T$ that starts at $v$, visits the neighbours of $w_1, w_2,\ldots, w_k$
on the right of $\cP$, and finishes at $v'$. The matching path of $v,v'$
is denoted by $P(v,v')$. Let $\cP'$ be the path obtained
from $\cP$ when substituting $[v,v']$ by $P(v,v')$. As shown in Figure~\ref{fig:upda_Path},
the path $\cP'$ goes from $W$ to $E$ and each vertex of $\cP'\backslash\{W,E\}$
has outdegree 2 on each side of $\cP'$ for the $\alpha_0$-orientation $O$.
Hence the path $\cP'$ is simple according to Lemma~\ref{lem:simple}. 
Moreover, by definition of $P(v,v')$, all edges of $T$
in the region enclosed by $[v,v']$ and $P(v,v')$ connect a vertex of 
$[v,v']\backslash \{v,v'\}$ to a vertex of $P(v,v')\backslash\{v,v'\}$,
see Figure~\ref{fig:upda}.

We can now describe the operations performed at each step of the sweeping 
process, as shown in Figure~\ref{fig:upda}.

\begin{itemize}
\item
Choose an admissible pair $v,v'$ of vertices on $\cP$.
\item
Color blue and orient from left to right all edges of $[v,v']$.
\item
Color red all edges inside the area enclosed by $[v,v']$ and $P(v,v')$,
and orient these edges from $P(v,v')$ to $[v,v']$.
\item
Update the path $\cP$, the part $[v,v']$ being replaced by $P(v,v')$.   
\end{itemize}
 According to the discussion above, the path $\cP$ is still simple after
these operations, and satisfies Invariant~(1). All the other invariants (2)-to-(6)
are easily shown to remain satisfied, as illustrated in Figure~\ref{fig:upda}.
At the end, the path $\cP$ is equal to $W\to S\to E$. The invariants~(3) and~(5)
ensure that the directions and colors of the inner edges of $T$ form
a transversal pair of bipolar orientations $Z$; 
and Invariant~(6) ensures that $Z$
is consistent with the $\alpha_0$-orientation $O$, i.e., $\Psi\circ\Phi(Z)=O$,
see
also Figure~\ref{fig:AlgoBij} for a complete execution of the algorithm.
This concludes the proof of Proposition~\ref{prop:lattice}.
 
\subsection{Essential circuits of an $\alpha_0$-orientation}
\label{sec:defTurnright}
Proposition~\ref{prop:lattice} 
ensures that the set $\cE$ of transversal edge-partitions of an irreducible 
triangulation is a distributive lattice, as $\cE$ is in bijection with the 
distributive lattice formed by the 
$\alpha_0$-orientations of the angular graph. 
By definition, the flip operation on $\cE$ is the effect of a 
flip operation on $\cO$ via the bijection. 
Recall that a flip operation on an $\alpha$-orientation consists in
reversing a clockwise essential circuit
(circuit with no chordal path). Hence, to describe 
the flip operation on $\cE$, we have to
characterise the essential circuits of an $\alpha_0$-orientation. 
For this purpose,
 we introduce 
the concept of \emph{straight path}. 

Consider an irreducible triangulation 
$T$  endowed with a transversal pair of bipolar orientations. 
Color blue the four outer edges of $T$ and orient them from $W$ to 
$E$. The conditions C1' and C2' ensure that there are four possible types 
for a bicolored angle $(e,e')$ 
of $T$, with $e'$ following $e$ in cw order: 
(outgoing red, outgoing blue) or (outgoing blue, ingoing 
red) or (ingoing red, ingoing blue) or (ingoing blue, outgoing red). The 
type of an edge of $Q(T)$ corresponding to a bicolored angle of $T$ 
(i.e., an edge going out of a black vertex) is 
defined as the type of the bicolored angle. For such an edge
$e$, the \emph{straight path} of $e$ is the 
oriented path $\mathcal{P}$ of edges of $Q(T)$ that starts at $e$ and such 
that each edge of $\mathcal{P}$ going out of a black vertex has the same type 
as $e$. Such a path is unique, as there is a unique choice for the 
outgoing edge at a white vertex (each white vertex has outdegree 1).

\begin{lemma}
\label{lemma:straight}
The straight path $\mathcal{P}$ of an edge $e\in Q(T)$ going out of a black 
vertex is simple and ends at an outer black vertex of $Q(T)$.
\end{lemma}
\begin{proof}
Notice that the conditions of a transversal pair of bipolar 
orientations remain satisfied if the directions of the edges of one color are 
reversed and then the colors of all inner edges are switched; hence the edge
 $e$ can be 
assumed to have type (outgoing red, outgoing blue) without loss of generality. 
Let $v_0,v_1,v_2,\ldots,v_i,\ldots$ be 
the sequence of vertices of the straight path $\mathcal{P}$ of $e$, so that 
the even indices correspond to black vertices of $Q(T)$ and the odd indices 
correspond to white vertices of $Q(T)$. Observe that, for $k\geq 0$, the 
black vertices $v_{2k}$ and $v_{2k+2}$ are adjacent in $T$ and the edge 
$(v_{2k},v_{2k+2})$ is either outgoing red or outgoing blue. Hence 
$(v_0,v_2,v_4,\ldots,v_{2k},\ldots)$ is an oriented path of $T$, so that it
is simple, according to Proposition~\ref{prop:acyc}. 
Hence, $\mathcal{P}$ does not pass twice by the same black vertex, i.e., 
$v_{2k}\neq v_{2k+2}$ for $k\neq k'$; and $\cP$ neither passes twice by the 
same white vertex (otherwise $v_{2k+1}=v_{2k+3}$ for 
$k\neq k'$ would imply $v_{2k+2}=v_{2k+4}$ by unicity of 
the outgoing edge at each white vertex, a contradiction). 
Thus $\mathcal{P}$ is a simple path, so that it
 ends at a black vertex of $Q(T)$ 
having no outgoing edge of type (outgoing red, outgoing blue), i.e., 
$\mathcal{P}$ ends at an outer black vertex of $Q(T)$. 
\hfill\phantom{1}  
\end{proof}

\begin{figure}
\begin{center}
\includegraphics[width=12cm]{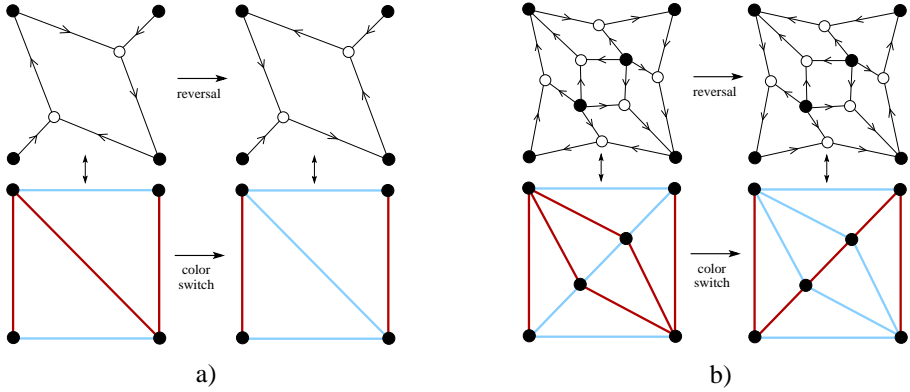}
\end{center}
\caption{The two possible configurations of an essential clockwise circuit 
$\mathcal{C}$ of $Q(T)$. In each case, an alternating 4-cycle is associated 
to the circuit; and reversing the circuit orientation 
corresponds to switching the edge
colors inside the alternating 4-cycle.}
\label{fig:cycle_reverse}
\end{figure}

\begin{proposition}
\label{prop:minimalcycle}
Given $T$ an irreducible triangulation and $Q(T)$ its angular graph
endowed with an $\alpha_0$-orientation $X$, 
an essential clockwise circuit $\mathcal{C}$ of $X$ 
satisfies either of the two following configurations,  
\begin{itemize}
\item
The circuit $\mathcal{C}$ is the boundary of a (quadrangular) inner 
face of $Q(T)$, see Figure~\ref{fig:cycle_reverse}(a).
\item
The circuit $\mathcal{C}$ has length 8. The four black vertices of 
$\mathcal{C}$ have no 
outgoing edge inside $\mathcal{C}$. The four white vertices 
of $\mathcal{C}$ have
their unique incident edge not on $\mathcal{C}$ inside 
$\mathcal{C}$, see Figure~\ref{fig:cycle_reverse}(b) for an example.
\end{itemize}  
\end{proposition}
\begin{proof}
First we claim that no edge of $Q(T)$ inside $\mathcal{C}$ 
has its origin on $\cC$; indeed the straight path construction ensures that
such an edge could be extended to a chordal path of $\cC$, which is 
impossible.
We define 
$\nb$, $\nw$ and $e$ as the number of black vertices, 
white vertices and edges inside $\mathcal{C}$. We denote by $2k$ the 
number of vertices 
on $\mathcal{C}$, so that there are $k$ black and $k$ white vertices on 
$\mathcal{C}$. Euler's relation and the fact that all 
inner faces of $Q(T)$ are quadrangular ensure that $(i):\ e=2(\nb+\nw)+k-2$.
As each white vertex of $Q(T)$ has degree 3, a white vertex on 
$\mathcal{C}$ has a unique 
incident edge not on $\mathcal{C}$. Let $l$ be the number of white 
vertices such that this 
incident edge is inside $\mathcal{C}$ (notice that $l\leq k$). 
Counting the edges inside
 $\mathcal{C}$ according to their incident white vertex gives 
$(ii):\ e=3\nw +l$. 
Edges 
inside $\mathcal{C}$ can also be counted
according to their origin for the $\alpha_0$-orientation. As 
no edge inside $\cC$ has its origin on $\cC$,
we have $(iii):\ e=4\nb +\nw$. 
Taking $2(i)-(ii)-(iii)$ yields $l=2k-4$. 
As $k$ is a positive integer and $l$ is a 
nonnegative integer 
satisfying $l\leq k$, the only possible values for $l$ and $k$ 
are $\{ k=4,l=4\}$, $\{ k=3,l=2\}$, and $\{k=2,l=0\}$. It is easily seen 
that the case $\{ k=3, l=2\}$ would correspond to a separating 3-cycle. 
Hence, the only possible cases are $\{ k=2,l=0\}$ and $\{k=4,l=4\}$,
shown in Figure~\ref{fig:cycle_reverse}(a) and~\ref{fig:cycle_reverse}(b),
respectively.
The first case corresponds to a circuit of length 4, which has to be the 
boundary of a face, as the angular graph of an irreducible 
triangulation (more generally, of a 3-connected plane graph) 
is well known to have no filled 4-cycle.
\hfill\phantom{1}
\end{proof}

\subsection{Flip operation on transversal structures}
\label{sec:fli}
As we prove now, the essential circuits of $\alpha_0$-orientations 
correspond to specific
patterns on transversal structures, making it possible to have a simple 
geometric interpretation of the flip operation, formulated directly
on the transversal structure.

Given $T$ an 
irreducible triangulation endowed with a transversal 
edge-partition, 
we define an \emph{alternating 4-cycle} as a cycle 
$\mathcal{C}=(e_1,e_2,e_3,e_4)$ of 4 edges of $T$ 
that are color-alternating (i.e., two adjacent edges of $\mathcal{C}$ 
have different colors). Given a vertex $v$ on $\mathcal{C}$, 
we call left-edge 
(right-edge) of $v$ the 
edge of $\mathcal{C}$ starting from $v$ and having the exterior of 
$\mathcal{C}$ on its 
left (on its right, respectively). 

\begin{lemma}
An alternating 4-cycle $\mathcal{C}$ in a transversal structure 
satisfies either of the
two following configurations.
\begin{itemize}
\item
All edges inside $\mathcal{C}$ and incident 
to a vertex $v$ of $\mathcal{C}$ have the color of the left-edge of 
$v$. Then $\mathcal{C}$ is called a \emph{left alternating 4-cycle}
\item
All edges inside $\mathcal{C}$ and incident to a vertex $v$ of $\mathcal{C}$ 
have the color of 
the right-edge of $v$. Then $\mathcal{C}$ is called a \emph{right
alternating 4-cycle}.
\end{itemize} 
\end{lemma}
\begin{proof}
Let $k$ be the number of vertices inside 
$\mathcal{C}$. 
Condition C1 ensures that there are $4k$ bicolored angles 
incident to a vertex inside 
$\mathcal{C}$.
Moreover, Euler's relation ensures that there are $2k+2$ faces 
inside 
$\mathcal{C}$. Hence Lemma~\ref{lemma:bic} implies that there are $4k+4$ 
bicolored angles inside $\mathcal{C}$. As a consequence, there are 
four bicolored angles inside $\mathcal{C}$ that are
 incident to a vertex of 
$\mathcal{C}$. As $\mathcal{C}$ is alternating, each of the four 
vertices of $\mathcal{C}$ is incident to at least one bicolored angle.
Hence the pigeonhole principle implies that each vertex of $\mathcal{C}$
is incident to one bicolored angle inside $\mathcal{C}$.
  Moreover, each of the four inner faces $(f_1,f_2,f_3,f_4)$ inside 
$\mathcal{C}$ and incident to 
an edge of $\mathcal{C}$ has two bicolored angles, 
according to  Lemma~\ref{lemma:bic}. As such a face $f_i$ has two angles 
incident to vertices of $\mathcal{C}$, at least one bicolored angle of $f_i$ 
is incident to a vertex of $\mathcal{C}$. As there are four bicolored angles 
inside $\mathcal{C}$ incident to vertices of $\mathcal{C}$, the pigeonhole 
principle ensures that each face $f_1,f_2,f_3,f_4$ has exactly one bicolored 
angle incident to a vertex of $\mathcal{C}$. As each vertex of $\mathcal{C}$ 
is incident to one bicolored angle inside $\mathcal{C}$, these angles are in 
the same direction. If they start from $\mathcal{C}$ in clockwise 
(counterclockwise) direction, 
then $\mathcal{C}$ is a right alternating 4-cycle (left alternating 4-cycle,
respectively).\hfill\phantom{1}
\end{proof}

\begin{theorem}
\label{theo:lattice}
The set of 
transversal edge-partitions of a fixed irreducible triangulation 
is a distributive lattice. 
The flip operation consists in switching the edge-colors 
inside a right alternating 4-cycle, turning it into a left 
alternating 4-cycle. 
\end{theorem} 
\begin{proof}
The two possible configurations for an essential clockwise 
circuit of an $\alpha_0$-orientation are 
represented respectively in 
Figure~\ref{fig:cycle_reverse}(a) and Figure~\ref{fig:cycle_reverse}(b).
In each case, the essential circuit of the angular graph 
corresponds to a right alternating 4-cycle 
on $T$. Notice that 
a clockwise face of $Q(T)$ corresponds to the 
vertex-empty alternating 4-cycle, 
whereas essential circuits of length 8  correspond 
to all possible alternating 4-cycles with 
at least one vertex in their interior (Figure~\ref{fig:cycle_reverse}(a) 
gives an example). 
As shown in Figure~\ref{fig:cycle_reverse}, the effect of reversing a clockwise
essential circuit of the angular graph is clearly 
a color switch of the edges inside the associated
alternating 4-cycle.
\hfill\phantom{1}\end{proof}

\vspace{0.2cm}

\noindent{\bf Definition.} Given an irreducible triangulation $T$, 
the transversal structure of $T$ with no right alternating 4-cycle is called
\emph{minimal}, as it is at the bottom of the distributive lattice.

\section{Bijection with ternary trees}
\label{sec:bij}
This section focuses on the description of  a bijection between ternary trees
and irreducible dissections, where the minimal transversal structure (for the 
distributive lattice) plays a crucial role.
The mapping from ternary trees to irreducible triangulations 
relies on so-called \emph{closure
operations}, as introduced by G. Schaeffer in his PhD~\cite{S-these}. 

\begin{figure}[t]
\begin{center}
\includegraphics[width=14cm]{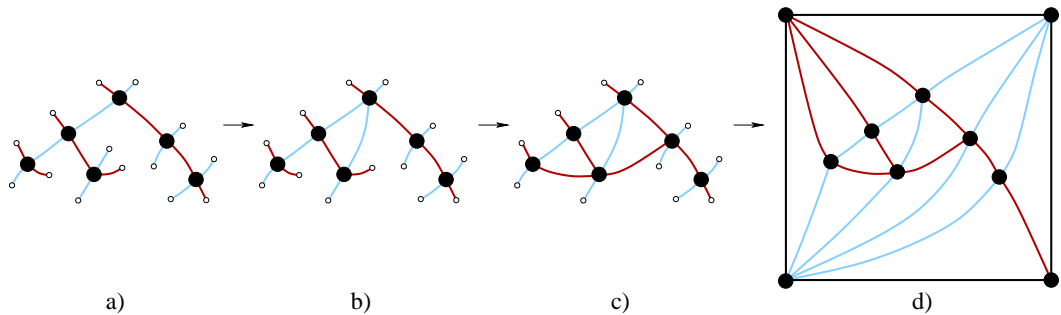}
\end{center}
\caption{The execution of the closure mapping on an example.}
\label{fig:Bij}
\end{figure}

\subsection{The closure mapping: from trees to triangulations}
\label{sec:closure}
A \emph{plane tree} is a plane graph with a unique face (the outer face). 
A ternary tree is a plane tree with vertex degrees in $\{1,4\}$. 
Vertices of degree 4 are called
\emph{nodes} and vertices of degree 1 are called \emph{leaves}. 
An edge of a ternary tree is called a \emph{closed edge} if it connects 
two nodes and
is called a \emph{stem} if it connects a node and a leaf.
It will be convenient to consider closed edges as made of two opposite 
half-edges meeting
at the middle of the edge, whereas stems will be considered as 
made of a unique half-edge incident to the  node 
and having not (yet) an opposite half-edge. A ternary tree
is \emph{rooted} 
by marking one leaf.
 The root allows us to distinguish the four neighbours of each  node, 
taken in ccw-order, into a parent 
(the neighbour in the direction of the root), 
a left-child, a middle-child, and a right-child. 
Thus, our definition of rooted 
ternary trees corresponds
 to the classical definition, where each node has three ordered children.

Starting from a ternary tree, the three steps to construct an 
irreducible triangulation are: 
local closure, partial closure and complete 
closure. Perform a 
counterclockwise walk alongside a ternary tree $A$ 
(imagine an ant walking around $A$ with the 
infinite face on 
its right). If a stem $s$ and then two sides of 
closed edges $e_1$ and $e_2$ are 
successively encountered 
during the traversal, create a half-edge opposite to the stem $s$ and
incident to the farthest extremity of $e_2$, 
so as to 
\emph{close} a triangular face. This operation is called a  
\emph{local closure}, see the transition Figure~\ref{fig:Bij}(a)-(b). 

The figure obtained in this way differs from $A$ by the presence of a 
triangular face and, 
more importantly, a stem $s$ of $A$ 
has become a closed edge, i.e., an edge made of two half-edges. 
Each time a sequence (stem, closed edge, closed edge) is found in ccw 
around the outer face
of the current figure, we perform a 
local closure, 
update the figure, and restart, until no local closure is possible. This 
greedy execution of 
local closures is called the \emph{partial closure} of $A$, see 
Figure~\ref{fig:Bij}(c). It is easily 
 shown that the figure $F$ obtained by performing the partial 
closure of $A$ does not 
depend on the order of 
execution of the local closures. Indeed, a
cyclic parenthesis word is associated to the counter-clockwise walk alongside
 the tree, with an opening parenthesis of weight 2 for a stem and
a closing parenthesis for a side of closed edge; the future local
closures correspond to matchings of the parenthesis word. 

At the end of the partial closure, the 
number $n_s$ of unmatched stems and the number $n_e$ of sides of closed edges 
incident to the outer face of $F$ satisfy the relation $n_s-n_e=4$. Indeed, 
this relation is satisfied on $A$ because a ternary tree with $n$  nodes 
has $n-1$ closed edges and $2n+2$ leaves (as proved by induction on the 
number of nodes); 
and the relation 
$n_s-n_e=4$ remains satisfied throughout the partial closure, as each local 
closure decreases $n_s$ and $n_e$ by 1. 
When no local closure is possible anymore, two consecutive 
unmatched stems on the boundary of the
outer face of $F$ are separated by at most one closed edge. 
Hence, the relation $n_s=n_e+4$ implies that the unmatched stems of $F$ are 
partitioned into four intervals $I_1,I_2,I_3,I_4$, where two consecutive 
stems of an interval are separated by one closed edge, and the 
last stem of $I_i$ is incident to the same vertex as the first stem of 
$I_{(i+1)\mathrm{mod}\ 4}$, see Figure~\ref{fig:Bij}(c). 

The last step of 
the closure mapping, called 
\emph{complete closure},
consists of the following operations. Draw a 4-gon $(v_1,v_2,v_3,v_4)$ 
outside of $F$; for $i\in\{1,2,3,4\}$, create an opposite half-edge for
 each  
stem $s$ of the interval $I_i$, the new half-edge being incident 
to the vertex $v_i$. Clearly, this process 
creates only triangular faces, so that it yields a triangulation 
of the 4-gon, see Figure~\ref{fig:Bij}(d).

Let us now explain how the closure mapping is related to
transversal edge-partitions. A ternary tree $A$ is said to be 
\emph{edge-bicolored} if each edge of $A$ (closed
edge or stem) is given a color ---red or blue--- such that any angle incident
to a  node of $A$ is bicolored, see 
Figure~\ref{fig:Bij}(a). 
Such a bicoloration, which is unique up to the choice of the colors, 
is called the 
\emph{edge-bicoloration} of $A$. 

\begin{lemma}
\label{lem:invariant}
Let $A$ be a ternary tree endowed with its edge-bicoloration.
The following invariant is maintained 
throughout the partial closure of $A$, 
\begin{center}
(I): \emph{any angle incident to the outer face of the current figure
is bicolored.}
\end{center}
\end{lemma}
\begin{proof}
By definition of the edge-bicoloration, $(I)$ 
is true on $A$. We claim that $(I)$ remains satisfied after each local 
closure. Indeed, let $(s,e_1,e_2)$ be the succession (stem, closed edge, 
closed edge)
 intervening in the local closure, let $v$ be the extremity of
 $e_2$ farthest from $s$, and let $e_3$ be the ccw 
follower of $e_2$ around $v$. Invariant $(I)$ implies that
 $s$ and $e_2$ have the same color. As we give to the new created 
half-edge $h$ the same color as its opposite half-edge $s$ 
(in order to have unicolored edges), $h$ and $e_2$ have the same color. 
The effect of the local closure on the angles of the outer face is the 
following: the angle $(e_1,e_2)$ disappears from 
the outer face, and the angle $(e_2,e_3)$ is replaced by the angle $(h,e_3)$. 
As $e_2$ has the same color as $h$, the bicolored angle $(e_2,e_3)$ is 
replaced by the bicolored angle $(h,e_3)$, so that $(I)$ remains true after 
the local closure.
\hfill\phantom{1}
\end{proof}

It also follows from this proof  that Condition C1
 remains satisfied throughout the partial closure, because the number of 
bicolored angles around each node is not increased, and is initially equal 
to 4. 
At the end of the partial closure, Invariant $(I)$ ensures that all stems of 
the intervals $I_1$ and $I_3$ are of one color, and all stems of the 
intervals $I_2$ and $I_4$ are of the other color. Hence, Condition C2 is 
satisfied after the complete closure, see Figure~\ref{fig:Bij}(d). Thus, 
the closure maps the edge-bicoloration of $A$ to a 
transversal edge-partition of the obtained triangulation of the 4-gon, 
which is in fact the 
minimal one: 

\begin{proposition}
\label{lemma:closure}
The closure of a ternary tree $A$ with $n$  nodes is an
irreducible triangulation $T$  with $n$ inner vertices.
The closure maps the edge-bicoloration of $A$ to the minimal 
transversal edge-partition of $T$.
\end{proposition}
\begin{proof}
Assume that $T$ has a separating 3-cycle $\mathcal{C}$. 
Observe that Lemma~\ref{lemma:bic} was stated and proved without the 
irreducibility condition. Hence, when the four outer edges of $T$ are 
colored blue, each inner face of $T$ has exactly two bicolored angles. 
Let $k\geq 1$ be the number of vertices inside $\mathcal{C}$. Euler's 
relation implies that the interior of $\cC$ contains 
$2k+1$ faces, so that 
there are $4k+2$ bicolored angles inside $\mathcal{C}$,
according to Lemma~\ref{lemma:bic}. Moreover, 
Condition C1 implies that there are 4k bicolored angles incident to a vertex 
that is in the interior of $\mathcal{C}$. 
Hence there are exactly two bicolored 
angles inside $\mathcal{C}$ incident to a vertex of $\mathcal{C}$. However, 
for each of the three edges $\{e_1,e_2,e_3\}$ of $\mathcal{C}$, the 
face incident to $e_i$ 
in the interior of $\mathcal{C}$ has at least one of its two bicolored 
angles incident to $e_i$. Hence, there are
at least 
three bicolored angles inside $\mathcal{C}$ and incident to a vertex of 
$\mathcal{C}$, a contradiction.

Now we show that the transversal edge-partition of $T$ induced by
the closure mapping is minimal, i.e., has no right alternating 4-cycle. Let 
$\mathcal{C}$ be an
alternating 4-cycle of $T$. This cycle has been closed during a local
closure involving one of the four edges of $\mathcal{C}$. Let $e$ be
this edge and let $v$ be the origin of the stem whose completion has
created the edge $e$. 
The fact
that the closure of a stem is always performed with the infinite face
on its right ensures that $e$ is the right-edge of $v$ on
$\mathcal{C}$, as defined in Section~\ref{sec:fli}. 
A second observation following from Invariant $(I)$ is that, when a stem $s$ 
is merged, the
angle formed by $s$ and by the edge following $s$ in
counterclockwise order around the origin of $s$ is a bicolored
angle. This ensures that $\mathcal{C}$ is a left alternating 4-cycle.
\hfill\phantom{1}\end{proof}

\subsection{Inverse mapping: the opening}
In this section, we describe the inverse of the closure mapping, 
from irreducible triangulations to ternary trees. As we have seen in the proof
of Lemma~\ref{lem:invariant},  during a local 
closure, the newly created half-edge $h$ has the same color as the 
clockwise-consecutive half-edge around the origin of $h$. 
Hence, for each half-edge $h$ incident to an inner vertex
of $T$,
\begin{itemize}
\item
if the angle formed by $h$ and its cw-consecutive half-edge is
unicolored, then $h$ has been created during a local closure,
\item
if the angle is bicolored, then $h$ is one of the 4
original half-edges of $A$ incident to $v$.
\end{itemize}

\begin{figure}
\begin{center}
\includegraphics[width=14cm]{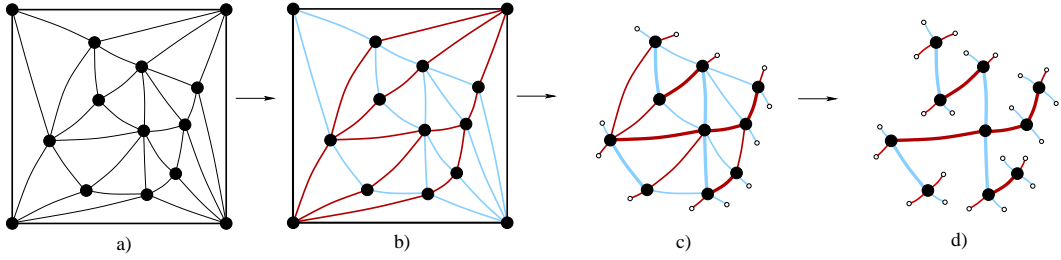}
\end{center}
\caption{The opening algorithm performed on an example.}
\label{figure:inverse}
\end{figure}

This property indicates how  to inverse  
the closure mapping. Given an irreducible triangulation $T$, the 
\emph{opening} of $T$ consists of the following steps, 
illustrated in Figure~\ref{figure:inverse}.

\begin{enumerate}
\item
Endow $T$ with its minimal transversal edge-partition.
\item
Remove the outer quadrangle of $T$ and all half-edges of $T$ incident to 
a vertex of the quadrangle.
\item
Remove all the half-edges whose clockwise-consecutive half-edge has the same color.
\end{enumerate}

The following lemma is a direct consequence of the definition of the 
opening mapping:
\begin{lemma}
\label{lemma:closure_inverse_opening}
Let $A$ be a ternary tree and let $T$ be the irreducible triangulation
 obtained by performing the closure of $A$. Then the opening of $T$
is $A$.
\end{lemma}

Hence, the closure $\Phi$ and the opening $\Psi$ are such that 
$\Psi\circ\Phi=\mathrm{Id}$. To prove that the opening and the closure 
mapping are mutually
 inverse, it remains to prove that $\Phi\circ\Psi=\mathrm{Id}$, which is 
more difficult and is done in two steps:\\ 
1) show that the opening of an irreducible triangulation $T$  
is a ternary tree,\\ 2) show 
that the closure of this ternary tree is $T$. 

The first step is to define an orientation of the half-edges of $T$ 
induced by the minimal transversal edge-partition. 
Each half-edge
of $T$ is associated to the angle on its right (looking from the incident 
vertex). We orient the half-edges of $T$ toward (outward of) their incident
vertex if the associated angle is unicolored (bicolored, respectively); 
with the restriction that half-edges on the outer 
quadrangle are leaved unoriented
and all angles incident to an outer vertex are considered as
unicolored.
This yields an orientation of the half-edges of $T$, which is called the 
\emph{4-orientation} of $T$.
Each inner vertex of $T$ is incident to four bicolored angle, hence
 has outdegree 4 in the 4-orientation. 
By definition of the 4-orientation, 
the opening of an irreducible  triangulation consists in removing the 
outer 4-gon and all ingoing
half-edges.

Let $e$ be an inner edge of $T$. An important remark is that the two 
half-edges  of $e$ can not be simultaneously
directed toward their respective incident vertex 
(otherwise, 
 the 4-cycle $\mathcal{C}$ bordering
the two triangular faces incident to $e$ would be a right
alternating 4-cycle, a contradiction).

Hence, only two cases arise for an inner edge $e$ of $T$.
\begin{itemize}
\item
If one half-edge of $e$ is ingoing, 
then $e$ is
called a \emph{stem-edge}. A stem-edge can be considered as simply oriented
(both half-edges have the same direction) for the 4-orientation.
\item
If the two half-edges of $e$ are 
 outgoing, then $e$ is called a
\emph{tree-edge}. A tree-edge can be considered as a bi-oriented edge for the
4-orientation. 
\end{itemize}

We define a clockwise circuit of the 4-orientation of $T$ as a simple cycle 
$\mathcal{C}$ of 
inner edges of $T$ such that each edge of $\mathcal{C}$ is either a 
tree-edge 
(i.e., a bi-oriented edge) or a stem-edge having the interior of $\mathcal{C}$ 
on its right.
\begin{lemma}
\label{lemma:noclockwise4ori}
The 4-orientation of $T$ has no clockwise circuit.
\end{lemma}
\begin{proof}
Assume there exists a clockwise circuit $\mathcal{C}$ in the 4-orientation 
of $T$. 
For a vertex $v$ on $\mathcal{C}$, we denote by $h_v$ the half-edge of 
$\mathcal{C}$ 
going out of $v$ when doing a clockwise traversal of $\mathcal{C}$, and
by $e_v$ the edge of $Q(T)$ following $h_v$ in clockwise order around $v$
(notice that $e_v$ is the most 
counterclockwise edge of $Q(T)$ incident to $v$ inside 
$\cC$). 
As $\mathcal{C}$ is a clockwise 
circuit for the 4-orientation of $T$, $h_v$ is directed 
outward of $v$. Hence the angle $\theta$ on the right of $h_v$ 
is a bicolored angle for the minimal transversal 
edge-partition of $T$, so that $e_v$ is going out of $v$ for the minimal 
$\alpha_0$-orientation $O_{\mathrm{min}}$ of $Q(T)$.

\begin{figure}
\begin{center}
\includegraphics{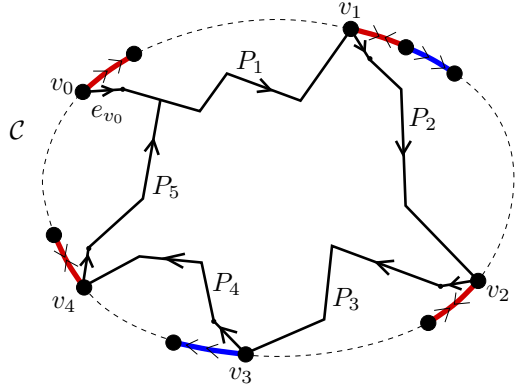}
\end{center}\caption{The existence of a clockwise circuit in the 4-orientation of $T$ 
implies the existence 
of a clockwise circuit in the minimal $\alpha_0$-orientation of $Q(T)$.}
\label{figure:cyclepaths}
\end{figure} 

We use this observation to build iteratively a clockwise circuit of 
$O_{\mathrm{min}}$ (see Figure~\ref{figure:cyclepaths}), 
yielding a contradiction. Let $v_0$ be a vertex on $\mathcal{C}$ and 
 $\mathcal{P}(v_0)$ 
the straight path starting at $e_{v_0}$, as defined in 
Section~\ref{sec:defTurnright}, for the minimal 
$\alpha_0$-orientation $O_{\mathrm{min}}$ of $Q(T)$. 
Lemma~\ref{lemma:straight} ensures 
that $\mathcal{P}(v_0)$ is simple and ends at an outer vertex of $Q(T)$.
In particular, $\mathcal{P}(v_0)$ has to 
reach $\mathcal{C}$ at a vertex $v_1$ different from $v_0$. 
We denote by $P_1$ the part of 
$\mathcal{P}(v_0)$ between $v_0$ and $v_1$, by $\Lambda_1$ 
the part of $\mathcal{C}$ between $v_1$ and $v_0$, and by $\mathcal{C}_1$ 
the cycle obtained by 
concatenating $P_1$ and $\Lambda_1$. Let $\mathcal{P}(v_1)$ be the 
straight path starting 
at $e_{v_1}$. The fact that $e_{v_1}$ is the most counterclockwise incident 
edge of $v_1$ in the 
interior of $\mathcal{C}$ ensures that $\mathcal{P}(v_1)$ starts in the 
interior of $\mathcal{C}_1$. 
The path $\mathcal{P}(v_1)$ has to reach $\mathcal{C}_1$ at a vertex 
$v_2\neq v_1$. 
We denote by $P_2$ the part of $\mathcal{P}(v_1)$ between 
$v_1$ and $v_2$. 
If $v_2$ belongs to $P_1$, then the concatenation of the part of $P_1$ 
between $v_2$ and $v_1$ and
of $P_2$ is a clockwise circuit of
$O_{\mathrm{min}}$, a 
contradiction. Hence $v_2$ is on $\Lambda_1$ strictly 
between $v_1$ and $v_0$. 
We denote by $\overline{P}_2$ 
the concatenation of $P_1$ and $P_2$, and by $\Lambda_2$ the part 
of $\mathcal{C}$ going from $v_2$
to $v_0$. As $v_2$ is strictly between $v_1$ and $v_0$, $\Lambda_2$ is 
strictly included in $\Lambda_1$. 
We denote by $\mathcal{C}_2$ the cycle made of the concatenation of 
$\overline{P}_2$ and $\Lambda_2$. Similarly as for the path 
$\mathcal{P}(v_1)$, 
the path $\mathcal{P}(v_2)$ must start in the interior of $\mathcal{C}_2$. 

Then we continue iteratively, see Figure~\ref{figure:cyclepaths}. At each 
step $k$, we consider the straight path $\mathcal{P}(v_k)$
 starting at $e_{v_k}$. This path starts in the interior of the cycle 
$\mathcal{C}_{k}$, and reaches again $\mathcal{C}_{k}$
 at a vertex $v_{k+1}$. The vertex $v_{k+1}$ can not belong to 
$\overline{P}_k$, 
otherwise a clockwise circuit of $O_{\mathrm{min}}$ would be created. 
Hence $v_{k+1}$ is
strictly between $v_{k}$ and $v_0$ on $\mathcal{C}$, i.e., is in 
$\Lambda_k\backslash\{ v_k,v_0\}$. 
In particular 
the path $\Lambda_{k+1}$ going from $v_{k+1}$ to $v_0$ on $\mathcal{C}$, is 
strictly included in the path $\Lambda_k$ going from $v_k$ to 
$v_0$ on $\mathcal{C}$. Thus, $\Lambda_k$ 
shrinks strictly at each step. Hence, there must be a step $k_0$ where 
$\mathcal{P}(v_{k_0})$ reaches $\mathcal{C}_{k_0}$ at a vertex on 
$\overline{P}_{k_0}$, thus creating a clockwise 
circuit of $O_{\mathrm{min}}$, a contradiction.\hfill\phantom{1}
\end{proof}    

\begin{lemma}
\label{lemma:spaning}
The tree-edges of $T$ form a tree spanning the inner vertices of $T$.
\end{lemma}
\begin{proof}
Denote by $H$ the graph consisting of the tree-edges of $T$ and their
incident vertices. A first observation is that $H$ has no cycle, as 
such a cycle of bi-oriented edges of $T$ would be a clockwise circuit in the 
4-orientation of $T$. 
Let $n$ be the number 
of inner vertices of $T$. Observe that $H$ can not be incident to the outer 
vertices of $T$, so that $H$ can cover at most the set of inner vertices of 
$T$.
 A well-known result of graph theory ensures that an acyclic 
graph $H$ having $n-1$ edges and
covering a subset of an $n$-vertex set $V$ is a tree covering 
exactly all vertices of $V$. Hence it remains to show that $H$ has $(n-1)$ 
edges. 
Let $s$ be the number of
stem-edges and $t$ be the number of tree-edges of $T$. As $T$ has $n$
inner vertices, there are $4n$ outgoing half-edges in the 4-orientation of 
$T$. 
Moreover, each stem-edge  has contribution 1 and each
tree-edge has contribution 2 to the number of outgoing half-edges. Hence,
$s+2t=4n$. Finally, Euler's relation ensures that $T$ has $(3n+1)$ inner
edges, so that $s+t=3n+1$. These two equalities ensure that $t=n-1$,
which concludes the proof that $H$ is a tree spanning the inner
vertices.
\hfill\phantom{1}\end{proof}

\begin{lemma}
The opening of an irreducible triangulation $T$  is a ternary tree.
\end{lemma}
\begin{proof}
As we have seen from the definition of the 4-orientation, 
the opening of an irreducible triangulation consists in 
removing the outer 4-gon and all ingoing
half-edges. The figure obtained in this way
 consists of the tree-edges, which form a 
spanning tree according to Lemma~\ref{lemma:spaning}, 
and of the half-edges that have lost their opposite half-edge. 
The edges of the first and second type correspond respectively
 to the closed edges and to the stems of the tree. 
In addition, after removing all ingoing half-edges, each vertex has degree 4, 
so that the tree satisfies the degree-conditions of a 
ternary tree.\hfill\phantom{1}  
\end{proof}

%To prove that the closure-mapping is a bijection whose inverse is the
%opening mapping, it remains to prove the following lemma:

\begin{lemma}
\label{lemma:opening_inverse_closure}
Let $T$ be an irreducible triangulation  and let $A$ be
the ternary tree obtained by doing the opening of $T$. Then
 the closure of $A$ is $T$.
\end{lemma}
\begin{proof}
First it is clear that
the complete closure (transition between 
Figure~\ref{fig:Bij}(c) and
Figure~\ref{fig:Bij}(d)) is the inverse of Step~2 of the opening algorithm. Let
$F$ be the figure obtained from $T$ after Step~2 of the opening
mapping. 

To prove that the partial closure of $A$ is $F$, it is sufficient
 to find a chronological order of deletion of the ingoing
half-edges of $F$ (for the 4-orientation) such that the inverse
of each half-edge deletion is a local closure. A local closure 
satisfies the property 
that the new created half-edge $h$ has the outer face on its 
right when $h$ is traversed toward its incident vertex.   
Thus the half-edge $h$
chosen to be disconnected at step $k$ must be incident to 
the outer face of the
current figure $F_k$, with $F_k$ on the
right when $h$ is traversed toward its incident vertex. We claim that
there always exists such a half-edge as long as 
there remain
 stem-edges in $F_k$. In that case, $F_k$
contains the spanning tree $H$ made of the tree-edges of $T$ plus at least
one stem-edge. Hence $F_k$ contains at least a cycle, thus a simple
cycle $\cC$ can be extracted from the boundary of $F_k$. There is at least one 
stem-edge $e$ on $\cC$, because no cycle is formed by tree-edges only. 
We claim that $e$ has the outer face of $F_k$ on its right, so that
the ingoing  half-edge $h$ of
$e$ is a candidate to be deleted (indeed if $e$
had the outer face of $F_k$ on its left, the concatenation
 of $e$ and of the path of tree-edges connecting the two extremities of $e$
 would be a clockwise circuit in the 4-orientation, a  
contradiction). As discussed above, the inverse operation
of the deletion
of $h$ is a local closure, which concludes the proof.
\hfill\phantom{1}\end{proof}

Finally, Lemma~\ref{lemma:closure_inverse_opening} and 
Lemma~\ref{lemma:opening_inverse_closure} yield the following theorem:

\begin{theorem}[bijection]
\label{theo:bij}
For $n\geq 1$, 
the closure mapping is a bijection between the set of ternary trees
with $n$  nodes and the set of irreducible triangulations with $n$ 
inner vertices. The inverse mapping of the closure is
the opening.
\end{theorem} 
 
\begin{theorem}[bijection, rooted version]
\label{cor:bije}
For $n\geq 1$,
the closure mapping induces a $(2n+2)$-to-4 correspondence between the set
$\mathcal{A}_n'$ of rooted ternary trees with $n$  nodes and the
set $\mathcal{T}_n'$ of rooted irreducible triangulations 
with $n$ inner vertices. In other words, $\mathcal{A}_n'\times \{1,\ldots,4\}$ 
is in 
bijection with $\mathcal{T}_n'\times \{1,\ldots,2n+2\}$.
\end{theorem}
\begin{proof}
It can easily be proved by induction on the number of  nodes that a 
ternary tree with $n$  nodes has $2n+2$ leaves. Hence, when rooting the 
ternary 
tree obtained by doing the opening of a triangulation in
 $\mathcal{T}_n'$, there are 
$2n+2$ 
possibilities to place the root. Conversely, starting from a rooted ternary 
tree with $n$  
nodes, there are four possibilities to place the root on the irreducible
triangulation  obtained by doing the closure of the tree, 
because the root has to be placed 
on one of 
the four outer edges.
\hfill\phantom{1}
\end{proof}

\subsection{Application to counting triangulations}
\label{sec:applic}
 In this section we focus on the application
to exact enumeration.  The key point
is that the bijection reduces the task of counting irreducible 
triangulations to the much easier
task of counting ternary trees. 
Then, we show that the enumeration of irreducible triangulations naturally leads to the 
enumeration of rooted 4-connected triangulations, which are closely related; 
the ingredients are
generating functions and a decomposition of a rooted irreducible triangulation as a sequence
of rooted 4-connected triangulations. 

\subsubsection{Counting irreducible triangulations}
As a first direct application, the bijection with ternary trees yields 
counting formulas for irreducible triangulations.

\begin{proposition}[counting irreducible triangulations]
\label{prop:count_irre}
For $n\geq 1$, the number of rooted irreducible triangulations 
with $n$ inner vertices is
$$
|\cT_n'|=4\frac{(3n)!}{n!(2n+2)!}.
$$
The number of unrooted irreducible triangulations 
with $n$ inner vertices is

\begin{eqnarray*}
|\cT_n|&=&\frac{(3n)!}{n!(2n+2)!}+\frac{1}{2}\frac{(3k)!}{k!(2k+1)!} \mathrm{\ \ \ \ \ \ \ \ \ \ \ \ \ \ \ \ \ \! if}\ n\equiv 0\!\!\mod 2\ \ [n=2k],\\[0.2cm]
|\cT_n|&=&\frac{(3n)!}{n!(2n+2)!}+\frac{1}{2}\frac{(3k+1)!}{k!(2k+2)!}+\frac{1}{2}\frac{(3k')!}{k'!(2k'+1)!}
\mathrm{\ \ \ \ \ \ \ \ if}\ n\equiv 1\!\!\mod 4\\
&&\ \ \ \ \ \ \ \ \ \ \ \ \ \ \ \ \ \ \ \ \ \ \ \ \ \ \ \ \ \ \ \ \ \ \ \ \ \ \ \ \ \ \ \ \ \ \ \ \ \ \ \ \ \ \ \ \  [n\!=\!2k\!+\!1\!=\!4k'\!+\!1],\\[0.2cm]
|\cT_n|&=&\frac{(3n)!}{n!(2n+2)!}+\frac{1}{2}\frac{(3k+1)!}{k!(2k+2)!} \mathrm{\ \ \ \ \ \ \ \ \ \ \ \ \ \! if}\ n\equiv 3\!\!\mod 4\ \ [n=1+2k].
\end{eqnarray*}
\end{proposition}
\begin{proof}
The enumerative formula follows from  
$|\mathcal{T}_n'|=\frac{4}{2n+2}|\mathcal{A}_n'|$ and 
from the well known fact that $|\mathcal{A}_n'|=(3n)!/((2n+1)!n!)$, 
which can be derived from the Lagrange inversion formula applied to the 
generating function $A(z)=z(1+A(z))^3$. The formula for $|\cT_n|=|\cA_n|$
follows from the enumeration of unrooted ternary trees, which is easily 
obtained by considering the possible rotation symmetries (order 2 around
a vertex or an edge, order 4 around a vertex). 
\hfill\phantom{1}
\end{proof}

The formula for \emph{rooted} irreducible triangulations can easily
be obtained from the series counting rooted 
triangulations of the 4-gon by using
a composition scheme, see~\cite{Tut}. 
To our knowledge, the counting formula for
\emph{unrooted} irreducible triangulations is new. 
However, a composition scheme should 
make it possible to count irreducible triangulations with a given
rotation symmetry (order 2 around a vertex or an edge, order 4 around
a vertex), starting from triangulations of the 4-gon with a given 
rotation symmetry, which have been counted by Brown~\cite{Br}.

\subsubsection{Counting rooted 4-connected triangulations}
\label{sec:4connected}
A graph is said to be \emph{4-connected}
 if it has more than three vertices and if at 
least four vertices 
have to be removed to disconnect it.
In this section we derive the enumeration of rooted 4-connected triangulations
from the counting formula for rooted irreducible triangulations. 
The idea is to translate 
a decomposition linking these two 
families of rooted triangulations to an equation linking their generating
functions. The net result we obtain is an explicit formula for 
the generating function
of rooted 4-connected triangulations (Proposition~\ref{prop:count}). 
We take advantage of the well known property that a triangulation 
is 4-connected 
iff the interior of any 3-cycle, except for the outer triangle, 
is a face.    In all this section, we denote by 
$\displaystyle\mathcal{T}'=\cup_n \mathcal{T}_n'$ and
by $\displaystyle\mathcal{C}'=\cup_n \mathcal{C}_n'$ the sets of rooted
irreducible triangulations and of rooted 
4-connected triangulations counted with respect to the number of 
inner vertices.

Observe the close connection between the 
definition of 4-connected 
triangulations and irreducible triangulations. 
In particular, for $n\geq 2$, the operation of 
removing the root
edge of an object of $\mathcal{C}_n'$ and carrying the root on the
counterclockwise-consecutive edge is an injective mapping from
$\mathcal{C}_n'$ to $\mathcal{T}_{n-1}'$. However, given $T \in
\mathcal{T}_{n-1}'$, the inverse edge-adding operation can create a
separating 3-cycle if there exists an internal path of length 2
connecting the origin of the root of $T$ to the vertex diametrically
opposed in the outer (quadrangular) face of $T$. Objects of $\mathcal{T}'$
having no such internal path are said to be \emph{undecomposable} and their
set, counted with respect to the number $n$ of inner vertices, is denoted by 
$\displaystyle\mathcal{U}'=\cup_n\mathcal{U}_{n}'$. 
The above discussion ensures
that $\mathcal{C}_n'$ is in bijection with $\mathcal{U}_{n-1}'$ for $n\geq 2$. 
In addition, 
a maximal
decomposition of an object $\gamma\in\mathcal{T}'$ along the above mentioned
interior paths of length 2 ensures that $\gamma$ is a sequence
of objects of $\mathcal{U}'$. Precisely, the graph enclosed by two 
consecutive paths of length 2 is either an undecomposable triangulation
 or is a quadrangle with a unique interior edge that connects the 
middles of the two paths. This leads to the equation
\begin{equation}
\label{eq:T}
T(z)+1=\frac{U(z)+1}{1-z(U(z)+1)},
\end{equation}
where $T(z)=\sum |\mathcal{T}_n'|z^n$ and $U(z)=\sum |\mathcal{U}_n'|z^n$ 
are respectively the series 
counting $\mathcal{T}'$ and $\mathcal{U}'$ with respect to the 
number of inner vertices.
\begin{proposition}
\label{prop:count}
The series $C(z)$ counting rooted 4-connected triangulations by their
number of inner vertices has the following expression,
\begin{equation}
\label{eq:C}
C(z)=\frac{z(A(z)-A(z)^2+1)}{1+z(A(z)-A(z)^2+1)},
\end{equation}
where $A(z)=z(1+A(z))^3$ is the series counting rooted ternary trees 
by their number of  nodes. 
\end{proposition}

\begin{proof}
As $\mathcal{C}_n'$ is in bijection with $\mathcal{U}_{n-1}'$ 
for $n\geq 2$ and 
as the unique 4-connected triangulation with less than 5 vertices is the 
tetrahedron, we have $C(z)=z(U(z)+1)$. Hence,
Equation~(\ref{eq:T}) yields $C(z)=z(T(z)+1)/(1+z(T(z)+1))$. Thus, it remains 
to provide an expression of $T(z)$ in terms of the series $A(z)=\sum A_nz^n$ 
counting rooted ternary trees by their number of  nodes.
 We define respectively the sets 
$\overline{\mathcal{A}_n}$ and $\widehat{\mathcal{A}_n}$ of ternary
trees with $n$
 nodes and having the following marks: a closed edge is marked and 
oriented for objects of
$\overline{\mathcal{A}}_n$, a closed edge is marked and oriented and a leaf is
marked for objects of $\widehat{\mathcal{A}}_n$. As a ternary tree
with $n$  nodes has $n-1$ closed edges and $2n+2$ leaves, we have
$A_n\cdot 2(n-1)=|\widehat{\mathcal{A}_n}|=|\overline{\mathcal{A}_n}|(2n+2)$, 
so that
$|\overline{\mathcal{A}_n}|=2\frac{n-1}{2n+2}A_n$. In addition, given a 
ternary tree with a marked oriented edge, the operation of cutting the 
marked edge produces an ordered pair of rooted ternary trees. Hence, the 
series counting
$\overline{\mathcal{A}}:=\cup_n \overline{\mathcal{A}}_n$ 
 with respect to the number of  nodes is
$A(z)^2$.
Finally, we know from Corollary~\ref{cor:bije} that
$T_n'=\frac{4}{2n+2}A_n$. Hence, we also have
$|\mathcal{T}_n'|=\frac{2n+2}{2n+2}A_n-2\frac{n-1}{2n+2}A_n=A_n-|\overline{\mathcal{A}_n}|$, from which we
conclude that $T(z)=A(z)-A(z)^2$. Finally, to obtain the expression of $C(z)$
in terms of $A(z)$, 
 substitute $T(z)$ by $A(z)-A(z)^2$ in the expression 
$C(z)=z(T(z)+1)/(1+z(T(z)+1))$. \hfill\phantom{1}\end{proof}

From Expression~(\ref{eq:C}), the coefficients of $C(z)$ 
can be quickly extracted: 
$C(z)=z+{z}^{3}+3\,{z}^{4}+12\,{z}^{5}+52\,{z}^{6}+241\,{z}^{7}+\mathcal{O}(z^8)$. 
The first coefficients match with the values given by the formula 
$\displaystyle c_n=\frac{1}{n}\sum_{i=1}^n
(-1)^i\binom{3n-i-1}{n-i}\binom{2i}{2}$
 for the number of rooted 4-connected triangulations with $n-1$ inner vertices. 
This formula was found by Tutte~\cite{Tut} using more complicated algebraic 
methods.

\section{Application to straight-line drawing}
\label{sec:draw}
A \emph{straight-line drawing} of a plane graph $G$ is 
a planar embedding of $G$ such that vertices are 
on an integer grid $[0,W]\times [0,H]$ and the edges are represented as segments.
The integers $W$ and $H$ are 
called the width and the height of the grid. As described next, 
transversal structures give rise to a simple 
straight-line drawing algorithm, based on face-counting operations.

\subsection{Face-counting algorithm}
Given an irreducible triangulation $T$ endowed with a 
transversal pair of bipolar orientations, we define the 
\emph{red-map} of $T$ as the plane graph $T_r$ obtained from $T$ 
by removing all
blue edges, see Figure~\ref{fig:Example}(e). The four outer edges of $T_r$
are additionally 
oriented from $S$ to $N$, so that $T_r$ is endowed with a bipolar
orientation, with $S$ and $N$ as poles (according to Corollary~\ref{prop:bip}).
Similarly, define the \emph{blue-map} $T_b$ of $T$ as the map obtained 
 by removing the red edges of $T$. The four outer edges of $T_b$ are 
additionally oriented from $W$ to $E$, so that $T_b$ is endowed with a  
 bipolar orientation. 
Given a vertex $v$ of $T$, we define the
\emph{leftmost outgoing red path} of $v$ as the oriented path $\prout(v)$ 
of edges
of the red-map starting at
$v$ and such that each edge of the path is the leftmost outgoing
edge at its origin.
As the orientation of the red-map is
bipolar, $\prout(v)$ is simple and ends at $N$. We also define the
\emph{rightmost ingoing red path} of $v$ as the path $\prin(v)$ 
starting from $v$ and
such that each edge of $\prin(v)$ is the rightmost ingoing edge at its 
end-vertex. 
This path is also acyclic and ends at $S$. The
\emph{separating red path} of $v$, denoted by $\cP_r(v)$, 
is the concatenation of $\prout(v)$ and $\prin(v)$. 
The oriented (simple) path $\mathcal{P}_r(v)$ goes from $S$ 
to $N$, thus inducing a bipartition of the inner faces of $T_r$, 
whether a face 
is on the left or on the right of $\mathcal{P}_r(v)$. Similarly, 
we define the 
\emph{leftmost outgoing blue path} $\pbout(v)$, 
the \emph{rightmost ingoing blue path} $\pbin(v)$, 
and write $\mathcal{P}_b(v)$ 
for the concatenation of $\pbout(v)$ and $\pbin(v)$, 
called the \emph{separating blue path} of $v$.

The 
following facts easily follow from the definition of 
the paths.

\begin{fact}
\label{fact:paths}
Let $v$ and $v'$ be two different vertices of $T$. 
%Let $\mathcal{P}_r^{\mathrm{out}}(v)$,  $\mathcal{P}_r^{\mathrm{out}}(v')$, 
%$\mathcal{P}_r^{\mathrm{in}}(v)$ and $\mathcal{P}_r^{\mathrm{in}}(v')$ 
%be respectively the leftmost outgoing red paths and rightmost ingoing red 
%paths of $v$ and $v'$. 
\begin{itemize}
\item
The paths  $\mathcal{P}_r^{\mathrm{out}}(v)$ and  
$\mathcal{P}_r^{\mathrm{out}}(v')$ do not cross each other, they join at a 
vertex $v''$ and then are 
equal between $v''$ and $N$.
\item
The paths  $\mathcal{P}_r^{\mathrm{in}}(v)$ and  
$\mathcal{P}_r^{\mathrm{in}}(v')$ do not cross each other, they join at a 
vertex $v''$ and then are 
equal between $v''$ and $S$.
\item
The paths $\mathcal{P}_r^{\mathrm{out}}(v)$ and 
$\mathcal{P}_r^{\mathrm{in}}(v')$ do not cross each other. 
In addition, $\mathcal{P}_r^{\mathrm{in}}(v')$ can not join 
$\mathcal{P}_r^{\mathrm{out}}(v)$ from the left of 
$\mathcal{P}_r^{\mathrm{out}}(v)$.  
\end{itemize}
\end{fact}

These facts easily imply that the separating 
red paths of two distinct vertices do not cross each other, and the same 
holds for separating blue paths.

The face-counting algorithm \textsc{TransversalDraw} runs as 
follows ---see Figure~\ref{fig:Triangulati3}---, 
the input being the
combinatorial description of an irreducible 
triangulation endowed with a transversal structure:

\begin{figure}
\begin{center}
\includegraphics[width=14cm]{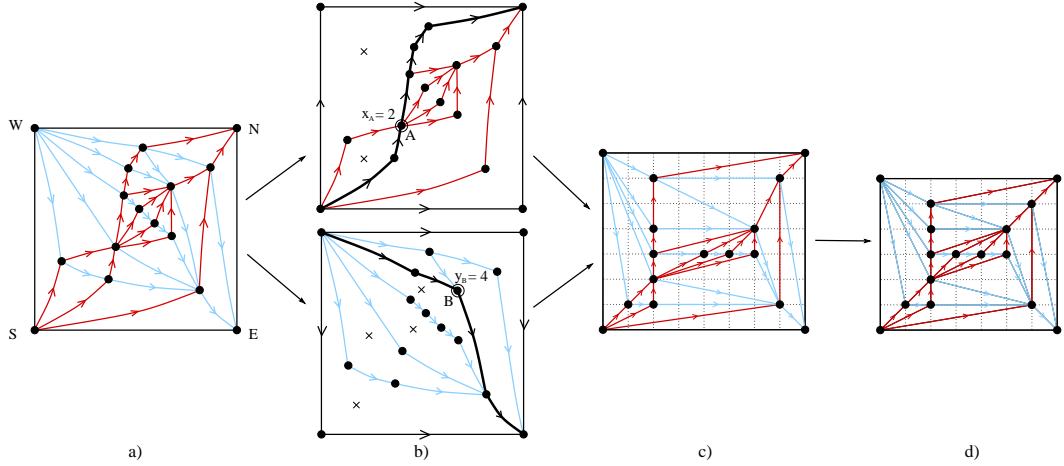}
\end{center}
\caption{The execution of \textsc{TransversalDraw} (Fig.a-to-c) 
on an example. 
The deletion of unused coordinates 
---\textsc{CompactTransversalDraw}--- is shown in Fig.d.}
\label{fig:Triangulati3}
\end{figure}

\vspace{0.6cm}

\begin{minipage}{14cm}
\noindent\textsc{TransversalDraw}($T$):

%\vspace{-0.1cm}

\begin{itemize}
%\item
%Place the outer vertices $S$, $W$, $E$, $N$ at 
%coordinates 
%$(0,0)$, $(0,f_b)$, $(f_r,0)$ and $(f_r,f_b)$, where $f_r$ and $f_b$ are
%the number of inner faces of the red-map and of the blue-map, respectively.
\item
For each vertex $v$ of $T$, %place $v$ on the grid in the following way:
\begin{itemize}
\item
the abscissa of $v$ is the number of inner faces of the red-map $T_r$ 
on the left of 
$\mathcal{P}_r(v)$,
\item
the ordinate of $v$ is the number of inner faces of the blue-map $T_b$ 
on the right of 
$\mathcal{P}_b(v)$.
\end{itemize}
\item
Place the vertices of $T$ on an integer grid accordingly, and link adjacent
vertices by segments.
\end{itemize}
\end{minipage}

\begin{theorem}
\label{theo:draw}
Given $T$ an irreducible triangulation endowed 
with a transversal structure,
\textsc{TransversalDraw}($T$) outputs  a straight-line drawing 
of $T$. The half-perimeter of the grid is $n-1$ if $T$ has $n$ vertices. 
The embedded edges satisfy the
following \emph{orientation property}:
\begin{itemize}
\item
red edges are oriented from bottom to top and weakly oriented from left 
to right,
\item
blue edges are oriented from left to right and weakly oriented from top 
to bottom.
\end{itemize}
The algorithm \textsc{TransversalDraw} can be implemented to run in 
linear time.
\end{theorem}
\begin{proof}
Let us first prove the orientation property.
Let $e=(v,v')$ be a red edge of $T$, directed from $v$ to $v'$. 
As discussed above, $\cP_b(v)$ and $\cP_b(v')$ do not cross each other.
Clearly Condition C1' implies that $\mathcal{P}_b(v)$ is on the
right of $\mathcal{P}_b(v')$. 
Hence, the ordinate of $v'$ is greater than the ordinate of $v$.
Then, proving that the abscissa of $v'$ is at least as large as 
the abscissa of $v$ reduces to proving that $\mathcal{P}_r(v)$ is not on
the right of $\mathcal{P}_r(v')$. This follows from two
easy observations: $\prout(v)$ is (weakly)
on the left of $\prout(v')$; and $\prin(v')$ is (weakly) on the right of 
$\prin(v')$.
Similarly, the blue edges are oriented from left to right and weakly
oriented from top to bottom. 

The proof that the embedding is a straight-line drawing relies 
on the orientation property and on the 
fact that the red and the blue edges 
are combinatorially transversal. To carry out the proof, we use
 a sweeping process akin to the 
iterative algorithm used to find the preimage of an $\alpha_0$-orientation
in Section~\ref{sec:corresTransAlpha}. 
The idea consists in maintaining an oriented 
blue path $\mathcal{P}$ from $W$ to $E$ called the \emph{sweeping path}, 
such that  the following invariant is maintained, 
$(I_{\mathrm{draw}})$:``the embedding of the 
edges of $T$ that are not topologically on the right of $\mathcal{P}$ is a straight-line 
drawing delimited to the top by the embedding of $(W,N)$, to the right 
by the embedding of $(N,E)$, and to the bottom-left by the embedding of 
$\mathcal{P}$''. The sweeping path is initially  $(W,N,E)$ and then 
``moves'' toward $S$. At each step, an inner face $f$ of the 
blue-map of $T$ is chosen, such that the left lateral path of $f$ 
is included in $\mathcal{P}$. 
Then, $\mathcal{P}$ is updated by replacing the left lateral path of $f$ by 
the right lateral path of $f$. Thus $\mathcal{P}$ remains an oriented path 
from $W$ to $E$ and is moved toward the bottom-left corner of the 
embedding, i.e., the vertex $S$. The fact that the invariants are 
maintained is easily checked from the orientation property
 and the fact that all the red edges inside $f$ 
connect transversally the two lateral blue paths, 
see Figure~\ref{fig:face_update}. At the end, the sweeping 
path is equal to $(W,S,E)$, so that the invariant $(I_{\mathrm{draw}})$
exactly implies that the embedding of $T$ is planar.

\begin{figure}
\begin{center}
 \includegraphics[width=12cm]{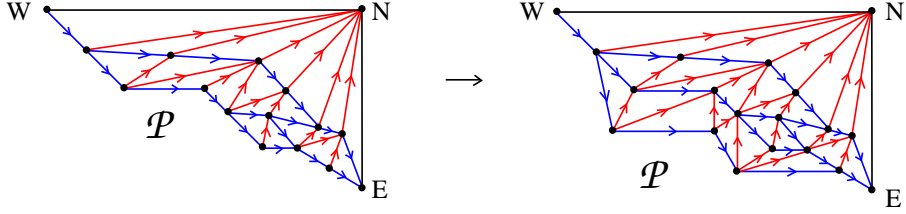}
\end{center}
\caption{The direction properties of the edges ensure that the drawing
is planar, step after step.}
\label{fig:face_update}
\end{figure}

We show now that the half-perimeter of the grid is equal to $n-1$
if $T$ has $n$ vertices. By definition of 
\textsc{TransversalDraw}, the minimal abscissa is $0$ and 
the maximal abscissa is equal to the number of inner faces $f_r$ of
the red-map. Similarly, the minimal ordinate is $0$ and the maximal 
ordinate is equal to the number of inner faces $f_b$ of the blue-map.  
Hence, the half-perimeter is equal to $f_r+f_b$. 
We write respectively $e_r$ and $e_b$ for
the number of edges of $T_r$ and $T_b$. Euler's relation ensures that
the total number $e$ of edges of $T$ is $3n-7$. Hence,
$e_b+e_r=e+4=3n-3$. In addition, Euler's relation, applied 
respectively to $T_r$
and $T_b$, ensures that $n+f_r=e_r+1$ and $n+f_b=e_b+1$. Thus,
$f_r+f_b=e_r+e_b-2n+2=n-1$. 

Finally, a linear implementation is obtained by suitably performing the 
face-counting operations.
For each inner vertex $v$ of $T$, consider the rightmost outgoing red path,
 the leftmost outgoing red path, the leftmost ingoing red path, and the
 rightmost ingoing red path. These four paths $P_1, P_2, P_3, P_4$ 
 partition the set of inner faces
 of the red-map $T_r$ into four areas 
$\mathcal{U}(v)$, $\mathcal{L}(v)$, $\mathcal{D}(v)$, and $\mathcal{S}(v)$,
 that are respectively enclosed by $(P_1,P_2)$, $(P_2,P_3)$, $(P_3,P_4)$, and $(P_4,P_1)$. 
Let $U(v)$, $L(v)$, $D(v)$ and $S(v)$ be the numbers of inner faces
 of $T_r$ in each of the four areas. The quantities $U(v)$ are
  easily computed in one pass, by  doing a traversal of the vertices
 of $T$ from $S$ to $N$. Similarly, the quantities $D(v)$ are
 computed in one pass doing a traversal from $N$ to $S$. 
Then, the quantities $L(v)$ are computed in one pass
 (using $D(v)$ and $U(v)$) by doing a traversal from $W$ to $E$. 
Finally, the abscissas of all vertices can be computed using 
$\mathrm{Abs}(v)=D(v)+L(v)$.  \hfill\phantom{1}\end{proof}

\subsection{Compaction step by coordinate-deletions}

Consider an irreducible triangulation $T$ endowed with a transversal structure 
and embedded using the face-counting algorithm
\textsc{TransversalDraw}. As observed in Figure~\ref{fig:Triangulati3}(c), 
some coordinate-lines might
be unoccupied. Hence, a natural further step is to delete the unused 
coordinates, yielding
a more compact drawing, as illustrated in Figure~\ref{fig:Triangulati3}(d). 
The corresponding algorithm is called 
\textsc{CompactTransversalDraw}.

\begin{theorem}
\label{theo:compactdraw}
Given $T$ an irreducible triangulation endowed 
with a transversal structure,
\textsc{CompactTransversalDraw}($T$) outputs  a straight-line drawing 
of $T$. The half-perimeter of the grid is at most $n-1$ if $T$ has 
$n$ vertices. The embedded edges satisfy the
same orientation property as \textsc{TransversalDraw}($T$), and the 
algorithm can be implemented
to run in linear time.
\end{theorem} 
\begin{proof}
Clearly, the orientation property of the red and of the blue edges remains 
satisfied after coordinate-deletions.
Notice that the planarity of \textsc{TransversalDraw} has been proved 
using only the orientation property, so 
that the same proof of planarity works as well for 
\textsc{CompactTransversalDraw}. The linear 
complexity follows from the easy property that coordinate-deletions 
can be performed in linear time.
\hfill\phantom{1} \end{proof}

\subsection{Drawing with the minimal transversal structure}
Notice that, in the definition of both \textsc{TransversalDraw} and 
\textsc{CompactTransversalDraw},
the irreducible triangulation is already equipped 
with a transversal structure. 
A complete straight-line drawing algorithm for irreducible 
triangulations is thus obtained
by computing a transversal structure first, and then 
launching the face-counting algorithm
\textsc{TransversalDraw}, possibly followed by the additional step
of deletion of unused 
coordinates (\textsc{CompactTransversalDraw}). 
A natural choice for the transversal structure computed is to take the
minimal one for the distributive lattice, with the further 
advantage of making it possible to perform
an average-case analysis of the grid size, as we will see 
in Section~\ref{sec:proofreduc}. The minimal 
transversal structure can be computed in linear time: 
1) compute a transversal structure in linear time
using an algorithm by Kant and He~\cite{Kant}, 
2) make the transversal structure minimal by 
iterated circuit reversals on the associated 
$\alpha_0$-orientations of the angular graph (the fact that
the overall complexity of circuit reversions is linear easily 
follows from ideas presented in~\cite{Khu93}).
Alternatively, a linear algorithm \textsc{ComputeMinimal} 
computing directly the minimal transversal structure is 
described in the PhD of the author~\cite{Fu-these}. 
We define the following straight-line drawing
algorithms for irreducible triangulations. 

\vspace{0.4cm}

\noindent\begin{minipage}{6cm}
\textsc{Draw}($T$):

%\vspace{-0.3cm}

\begin{itemize}
\item
Call \textsc{ComputeMinimal}($T$)
\item
Call \textsc{TransversalDraw}($T$)
\end{itemize}
\end{minipage}
\hspace{0.1cm}
\begin{minipage}{8cm}
\textsc{CompactDraw}($T$):

%\vspace{-0.3cm}

\begin{itemize}
\item
Call \textsc{ComputeMinimal}($T$)
\item
Call \textsc{CompactTransversalDraw}($T$)
\end{itemize}
\end{minipage}

\begin{figure}
\begin{center}
\includegraphics[width=10cm]{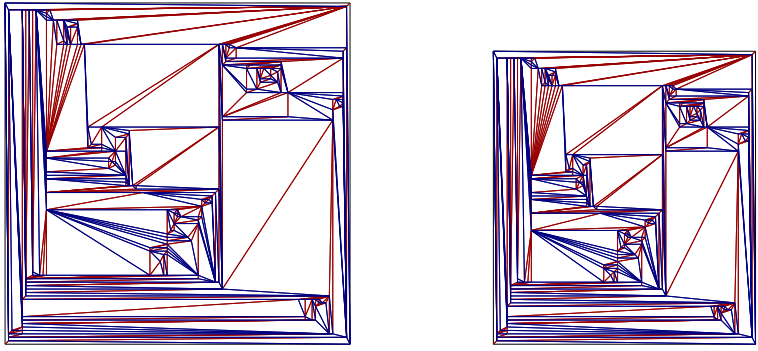}
\end{center}
\caption{A random triangulation with 200 vertices embedded with the algorithms
\textsc{Draw} 
and \textsc{CompactDraw}.}
\label{fig:big}
\end{figure}
%As we have seen, the closure mapping gives rise to an
%efficient (linear time) uniform sampler \textsc{SampleT}($n$) on 
%$\mathcal{T}_n$. In addition, sampled objects of $\mathcal{T}_n$
%are naturally endowed, through the closure, with their minimal
%transversal edge-partition. Hence, we can easily run the face-counting
%algorithms \textsc{TransversalDraw} and
%\textsc{CompactTransversalDraw} on the sampled objects. 

\vspace{0.2cm}

Simulations on random irreducible triangulations of 
large size ($n\approx 50000$)
indicate that the grid size is always
approximately $\frac{n}{2}\times \frac{n}{2}$ with
\textsc{Draw} and $\frac{n}{2}(1-\alpha)\times
\frac{n}{2}(1-\alpha)$ with \textsc{CompactTransversalDraw}, for some
constant $\alpha\approx 0.18$, see Figure~\ref{fig:big} for an example. 
It turns out that the  size of the grid can be
 analyzed thanks to the closure mapping presented in Section~\ref{sec:bij}. In the same way, 
as described by Bonichon \emph{et al}~\cite{BGHPS04}, 
the bijection of~\cite{PS03} makes it possible to perform 
a probabilistic analysis of the grid size for
straight-line drawing algorithms based on Schnyder woods~\cite{Fel}. 

Finally let us mention that a more general procedure of compaction is described in the PhD
of the author~\cite{Fu-these}: (1) delete a set $\cE$ of red and blue edges such that the red
acyclic orientation and the blue acyclic orientation remain bipolar, (2) call the face-counting
algorithm \textsc{TransversalDraw}, (3) re-insert the edges of $\cE$. It is shown in~\cite{Fu-these}
that the drawing thus obtained is planar and has semi-perimeter $n-1-|\cE|$. (However the 
orientation properties might not be satisfied.) 
%Indeed, the unused abscissas and ordinates of 
%\textsc{TransversalDraw} correspond 
%to particular inner edges of the ternary tree, whose number can be proven to 
%be 
%asymptotically concentrated around  $5n/27$ for a random ternary tree
%with $n$  nodes.

\section{Analysis of the grid size of the drawing algorithms}
\label{sec:proofreduc}
In this section, we explain how the closure mapping gives 
rise to a precise probabilistic 
analysis of the grid size of the algorithms 
\textsc{Draw} and \textsc{CompactDraw}. The main result is as follows.

\begin{theorem}
\label{prop:reduction}
Let $T$ be a rooted irreducible triangulation with $n$ inner 
vertices taken uniformly at random. 
Let $W\times H$ be the grid size of \textsc{Draw}($T$)
and let $W_c\times H_c$ be the grid size of 
\textsc{CompactDraw}($T$). Then, the following results hold 
asymptotically with high probability, 
up to fluctuations of order $\sqrt{n}$,
$$W\times H\simeq \frac{n}{2}\times \frac{n}{2},\ \ \ \ \ \  W_c\times H_c\simeq \frac{11}{27}n\times \frac{11}{27}n.$$ 

%The same result holds with an 
%$\epsilon$-formulation: for any fixed $\epsilon>0$, the probability that 
%$W$ or $H$ are outside of 
%$[\frac{n}{2}(1-\epsilon), \frac{n}{2}(1+\epsilon)]$ and the probability that 
%$W_c$ or $H_c$ are outside of $[\frac{11n}{27}(1-\epsilon), \frac{11n}{27}(1+\epsilon)]$ are 
%asymptotically exponentially small.
\end{theorem}
Let us mention that the concentration property  also holds in a so-called \emph{$\epsilon$-formulation}: ``for any $\epsilon>0$, the probability that $W$ or $H$ are outside of $[(1/2-\epsilon)n,(1/2+\epsilon)n]$ and the probability that  $W_c$ or $H_c$ are outside of $[(11/27-\epsilon)n,(11/27+\epsilon)n]$ are asymptotically exponentially small''. 

\vspace{0.5cm}

We first deal with the analysis of $W$ and $H$, i.e., the width and the 
height of 
\textsc{Draw}($T$). This task is rather easy and allows us to 
introduce some tools 
that will also be used to analyze the grid size of 
\textsc{CompactDraw}($T$), 
namely generating functions and the so-called quasi power theorem.  
In the sequel,
the set of rooted irreducible triangulations with $n$ inner  
vertices is denoted by $\cT_n'$.

\subsection{Analysis of the grid size of \textsc{Draw}($T$)}
Given $T\in\mathcal{T}_n'$ endowed with its minimal transversal pair of 
bipolar orientations, the 
width of the grid of \textsc{Draw}($T$) is, by definition, the 
number of inner faces 
of the red-map $T_r$. Let $e_r$ and $f_r$ be the numbers of inner edges and
of inner faces of $T_r$. Euler's relation applied to $T_r$ ensures 
that $e_r=n+f_r-1$. 
By definition 
of the opening, $e_r$ is equal to the number of red edges 
(including the stems) of the edge-bicolored ternary 
tree obtained by doing the opening of $T$. An easy adaptation of the proof of 
Theorem~\ref{cor:bije} ensures that the uniform distribution on 
$\mathcal{T}_n'$ is mapped by the opening to 
the uniform distribution on rooted edge-bicolored ternary tree with $n$ 
 nodes. These 
observations lead to the following statement:
\begin{fact}
\label{fact:distr}
The distribution of the width of \textsc{Draw($T$)} for $T$ 
uniformly sampled in 
$\mathcal{T}_n'$ is equal to the distribution of $e_r-n-1$, where $e_r$ is 
the number of red 
edges of a uniformly sampled rooted edge-bicolored ternary tree with $n$ 
 nodes.
\end{fact}

We denote by $\mathcal{R}$ (\emph{resp.} $\mathcal{B}$) the set of rooted 
edge-bicolored ternary trees 
whose root leaf is incident to a red stem (\emph{resp.} 
blue stem). 
For a rooted 
edge-bicolored ternary tree $\gamma$, we denote by $|\gamma|$ the number of  
nodes of $\gamma$ 
and by $\xi(\gamma)$ the number of red edges of $\gamma$. We define the 
generating functions 
$$R(z,u)=\sum_{\gamma\in\mathcal{R}}z^{|\gamma|}u^{\xi(\gamma)},\ \ B(z,u)=\sum_{\gamma\in\mathcal{B}}z^{|\gamma|}u^{\xi(\gamma)}$$ 
that count the set 
$\mathcal{R}$ and the set $\mathcal{B}$ with respect to the number of  
nodes and the number of red edges. The generating function $E(z,u)$ counting 
rooted edge-bicolored ternary trees with respect to the number of  nodes 
and the number of red edges is thus equal to $R(z,u)+B(z,u)$. The classical 
decomposition of a rooted ternary tree at the root node into three subtrees 
 translates to the following equation system:
$$
\left\{
\begin{array}{rcl}
R(z,u)&=&zu\left(1+B(z,u)\right)\left(u+R(z,u)\right)\left(1+B(z,u)\right),\\
B(z,u)&=&z\left(u+R(z,u)\right)\left(1+B(z,u)\right)\left(u+R(z,u)\right).
\end{array}\right.
$$
This system is polynomial in $R$, $B$, $z$ and $u$, so that one derives, 
 by algebraic elimination, a trivariate polynomial 
$P(E,z,u)$ that satisfies $P(E(z,u),z,u)=0$, i.e., $E(z,u)$ is an 
algebraic series.

We state now an adaptation of the quasi-power theorem~\cite[Thm. IX.7,
Cor. IX.1]{fla},~\cite{drmota97systems} 
for algebraic 
series.  This theorem makes it possible to analyse the distribution of the number of red edges in a 
rooted edge-bicolored ternary tree.

\begin{theorem}[Algebraic quasi power theorem] Let $E(z,u)$ be the generating 
function of a combinatorial class $\mathcal{E}$, where the variable $z$ marks 
the size $|\gamma|$ of an object $\gamma\in\mathcal{E}$ and the variable $u$ 
marks a parameter $\xi$, i.e., 
$E(z,u)=\sum_{\gamma\in\mathcal{E}}z^{|\gamma|}u^{\xi(\gamma)}$. Assume that 
$E(z,u)$ is an algebraic series, i.e., there exists a trivariate 
polynomial $P(E,z,u)$ with rational coefficients such that $P(E(z,u),z,u)=0$. 
Consider the polynomial system:
$$
\mathrm{(S)}:=\{  P(E,z,u)=0,\ \ \ \frac{\partial P}{\partial E}(E,z,u)=0\}.
$$
Assume that $(S)$ has a solution $\{\tau,\rho\}$ at $u=1$ such that $\tau$ and 
$\rho$ are positive real values and there exists no other complex 
solution $(E,z)$ of $(S)$ at $u=1$ such that $|z|\leq \rho$. Assume further 
that the \emph{derivative-condition} 
$\{\frac{\partial^2 P}{\partial E^2}(\tau,\rho,1)\neq 0,\frac{\partial P}{\partial z}(\tau,\rho,1)\neq 0\}$ is satisfied. Assume also that the Jacobian 
of $(S)$ with respect to $E$ and $z$ does not vanish at $(\tau,\rho,1)$, i.e.,
$$
\mathrm{det}\left( \begin{array}{cc}\frac{\partial P}{\partial E}(\tau,\rho,1)&\frac{\partial^2 P}{\partial E^2}(\tau,\rho,1)\\
\frac{\partial P}{\partial z}(\tau,\rho,1)&\frac{\partial^2 P}{\partial z\partial E}(\tau,\rho,1)
\end{array}\right)\neq 0.
$$
Then there exists a unique pair of algebraic series $(\tau(u)$, $\rho(u))$ 
such 
that $$\{P(\tau(u),\rho(u),u)=0,\frac{\partial P}{\partial E}(\tau(u),\rho(u),u)=0, \tau(1)=\tau,\rho(1)=\rho\}.$$ Finally, assume that $\rho(u)$ 
satisfies the 
mean condition $\rho'(1)\neq 0$ and variance condition 
$\rho''(1)\rho(1)+\rho'(1)\rho(1)-\rho'(1)^2\neq 0$. Then, for $\gamma$ taken 
uniformly at 
random in the set of objects 
of $\mathcal{E}$ of size $n$, the random variable 
$X_n=\xi(\gamma)$ is asymptotically equal to 
$\mu n$, up to Gaussian fluctuations of order $\sigma \sqrt{n}$, where 
$$\mu=-\frac{\rho'(1)}{\rho(1)},\ \ \sigma^2=-\frac{\rho''(1)}{\rho(1)}-\frac{\rho'(1)}{\rho(1)}+\left( \frac{\rho'(1)}{\rho(1)}\right)^2.$$ 
In other words, $(X_n-\mu n)/(\sigma\sqrt{n})$ converges 
to the normal law. 
%The result
%also holds with an $\epsilon$-formulation: for any $\epsilon >0$, the 
%probability that $X_n$ is outside of $[(1-\epsilon)\mu n,(1+\epsilon)\mu n]$ 
%is asymptotically exponentially small.
\end{theorem}
This theorem, despite a rather long statement and several conditions to 
check, is easy to apply in practice. In our case, $E(z,1)$ is clearly equal 
to $2A(z)$ with $A(z)=z(1+A(z))^3$ the generating function of rooted ternary 
trees. Hence, at $u=1$, we have $P(E,z)=E/2-z(1+E/2)^3$. Using a computer 
algebra software, the solution of $(S)$ at $u=1$ is found to be 
$\{ \tau=1, \rho=4/27\}$. The derivative condition and the Jacobian 
condition are then easily checked. The algebraic function $\rho(u)$ is 
obtained by taking the resultant of the two equations of $(S)$ with respect 
to $E$; then the factor $Q(z,u)$ of the resultant such that $Q(\rho,1)=0$ 
gives an algebraic equation for $\rho(u)$, i.e., $Q(\rho(u),u)=0$. 
From the algebraic equation $Q(\rho(u),u)=0$, the derivative and 
second derivative of $\rho(u)$ at $u=1$ are readily calculated. 
To calculate $\rho'(1)$, we derive $Q(\rho(u),u)=0$ and find 
$\frac{\partial Q}{\partial z}(\rho(u),u)\rho'(u)+\frac{\partial Q}{\partial u}(\rho(u),u)=0$, so that
$\rho'(1)=-\frac{\partial Q}{\partial u}(\rho,\tau)/\frac{\partial Q}{\partial z}(\rho,\tau)$. 
The 
mean condition and variance condition are readily checked and we find 
$\mu=3/2$. Hence, the number $e_r$ of red edges in a random rooted 
edge-bicolored ternary tree is asymptotically equal to $3n/2$, up to 
fluctuations of order $\sqrt{n}$.
This result, together with Fact~\ref{fact:distr}, ensure that the width of 
\textsc{Draw}($T$), for $T$ taken uniformly at random in $\mathcal{T}_n'$, 
is asymptotically equal to $n/2$ up to fluctuations of order 
$\sqrt{n}$.
The result is the same for the height $H$ of \textsc{Draw}($T$),
by stability of $\mathcal{T}_n'$ and of the algorithm \textsc{Draw} 
under the $\pi/2$-clockwise rotation.   Thus, we have proved 
the statement of Theorem~\ref{prop:reduction} on the distribution of the 
grid size of \textsc{Draw}($T$).

\subsection{Analysis of the grid size of \textsc{CompactDraw}($T$)}
\subsubsection{Overview of the proof}
Now we concentrate on the distribution of the grid width of 
\textsc{CompactDraw}($T$), for $T$ taken uniformly at random in 
$\mathcal{T}_n'$ (the analysis of the grid height is similar).  
By definition, the width of 
\textsc{CompactDraw}($T$) is $W-\Delta_{\mathrm{abs}(T)}$, 
where $W$ is the width of \textsc{Draw}($T$) 
and $\Delta_{\mathrm{abs}}$ is the number of abscissas not used 
by \textsc{Draw}($T$). We have already proved that $W$ 
is asymptotically equal to $n/2$ up to fluctuations of 
order $\sqrt{n}$.
Hence, to obtain the statement of Theorem~\ref{prop:reduction} 
about $W_c$, it is sufficient to prove that $\Delta_{\mathrm{abs}(T)}$ is 
asymptotically equal to $5n/54$ up to fluctuations of order 
$\sqrt{n}$.
The steps of the proof are the following: first, the abscissas that are not 
used by 
\textsc{Draw}($T$) are associated with 
specific edges of $T$; then these edges of $T$ are shown to 
correspond to specific red edges of the ternary tree obtained by doing the 
opening of $T$. Finally, using generating functions and the algebraic quasi 
power theorem, it is proved that the number of such edges in a random rooted 
edge-bicolored
ternary tree with $n$  nodes is asymptotically equal to 
$5n/54$, up to fluctuations of order $\sqrt{n}$.

\subsubsection{Characterisation of unused abscissas as particular edges 
of $T$}
We analyze the number of unused abscissas of \textsc{Draw}($T$), 
i.e., the number of vertical lines of the grid that bear no vertex of $T$. 
Recall that the width of the grid of \textsc{Draw}($T$) is the 
number of inner faces $f_r$ of the red-map $T_r$ of $T$, where $T_r$ is 
obtained 
by computing the minimal transversal pair of bipolar orientations of $T$ 
and then removing all blue edges. By definition of
\textsc{TransversalDraw}($T$),
the abscissa $\mathrm{Abs}(v)$ of a vertex $v$ of $T$ is obtained by 
associating with $v$ an oriented path $\mathcal{P}_r(v)$ of the red-map, called 
the separating red path of $v$ (as defined in Section~\ref{sec:draw}), and then 
counting 
the number of inner faces of $T_r$ on the left of $\mathcal{P}_r(v)$. 
We show first that each 
abscissa-candidate $1\leq i\leq f_r$ can be associated with a unique 
inner face of $T_r$ ---for which we write $\mathrm{Abs}(f)=i$---
such that the existence of a vertex $v$ of $T$ with $\mathrm{Abs}(v)=i$ 
only depends on the configuration of the edge at the bottom-right corner 
of $f$.

Let us start with a few definitions. Given $e=(v,v')$ an edge of 
$T_r$ directed from $v$ to $v'$, the \emph{separating red path} of $e$, 
denoted by $\mathcal{P}_r(e)$, 
is the path obtained by concatenating 
the leftmost 
outgoing red path of $v'$, the edge $e$, and the rightmost 
ingoing 
red path of $v$. We recall that, for each inner face $f$ of $T_r$, there
are two vertices $s_f$ and $t_f$ of $f$ such that
 the boundary of $f$ consists of a left and a right lateral paths going 
from $s_f$ to $t_f$.
The \emph{separating red path} of $f$, denoted by
$\mathcal{P}_r(f)$, is defined as the concatenation of the leftmost 
outgoing red path of $t_f$, of the 
right lateral path $P_{\mathrm{right}}(f)$ of $f$, 
and of the rightmost ingoing 
red path of $s_f$. 
Observe that the first edge $e_f$ 
of $P_{\mathrm{right}}(f)$, called bottom-right 
edge of $f$, satisfies $\mathcal{P}_r(f)=\mathcal{P}_r(e_f)$.

\begin{lemma}
\label{lemma:cross}
Let $e$ and $e'$ be two different edges of the red-map $T_r$. 
Then the separating paths 
$\mathcal{P}_r(e)$ and $\mathcal{P}_r(e')$ do not cross each other.
\end{lemma}
\begin{proof}
First it is easily checked that the bipolar orientation on 
red edges (and similarly the one
on blue edges) has 
no transitive edge; otherwise, there would be a face of the 
red-map with one  lateral path reduced
to an edge, in contradiction with the presence of transversal 
blue edges within the face.
 Hence only three cases can arise: either $e'$ and 
$\mathcal{P}_r(e)$ do not intersect, or they intersect at a unique extremity 
of $e'$, or $e'$ is on $\mathcal{P}_r(e)$. Fact~\ref{fact:paths} allows us to 
check that $\mathcal{P}_r(e)$ and $\mathcal{P}_r(e')$ do not cross each other 
for each of the 
three cases. \hfill\phantom{1}
\end{proof}

Recall that the separating red path of a face $f$ is the separating red path 
of its bottom-right edge. Thus, Lemma~\ref{lemma:cross} ensures that the 
separating 
 red paths $\mathcal{P}_r(f)$ and $\mathcal{P}_r(f')$ of two different inner 
faces 
$f$ and $f'$ of 
$T_r$ do not cross each other. In addition, it is easy to see that they are 
different using the fact that the bottom-right edge of a face is not the 
leftmost outgoing red edge at its origin. As a consequence, 
$\mathrm{Abs}(f)\neq \mathrm{Abs}(f')$. There are $f_r$ inner faces in $T_r$, 
each inner face $f$ clearly satisfying $1\leq \mathrm{Abs}(f)\leq f_r$. 
Hence the pigeonhole principle yields:

\begin{fact}
\label{fact:abscissa}
For $1\leq i\leq f_r$, there exists a unique inner face $f$ of $T$ such that 
$\mathrm{Abs}(f)=i$.
\end{fact}

Thus, the inner faces 
of $T_r$ are strictly ordered from left to right according to their 
associated abscissa.

\begin{lemma}
\label{lemma:sepvertedge}
The separating red paths of an edge $e$ of the red-map $T_r$ and of a
 vertex $v\in T$ do not cross each other. In addition, given 
$e$ an edge of $T_r$, there exists a vertex 
$v$ 
such that $\mathcal{P}_r(v)=\mathcal{P}_r(e)$ iff either $e$ is the 
rightmost ingoing red edge at its end-vertex or $e$ is the leftmost 
outgoing red edge at its origin. 
\end{lemma}
\begin{proof}
The fact that $\mathcal{P}_r(e)$ and $\mathcal{P}_r(v)$ do not cross each 
other is easily 
checked using 
Fact~\ref{fact:paths}. The second statement of the lemma follows from a few 
observations. Denote by $v_1$ the origin of $e$ and by $v_2$ the
end-vertex of $e$. If $v$ is not on $\mathcal{P}_r(e)$ then clearly 
$\mathcal{P}_r(v)$ is not equal to $\mathcal{P}_r(e)$. If $v$ is on 
$\mathcal{P}_r(e)$ between $v_2$ and $N$, then 
$\mathcal{P}_r(v)=\mathcal{P}_r(e)$ iff all edges of $\mathcal{P}_r(e)$ 
between $v_1$ and $v$ are the rightmost ingoing red edge at their end-vertex. 
If $v$ is on $\mathcal{P}_r(e)$ between $S$ and $v_1$, then 
$\mathcal{P}_r(v)=\mathcal{P}_r(e)$ iff all edges of $\mathcal{P}_r(e)$ 
between $v$ and $v_2$ are the leftmost outgoing red edge at their origin. 
It follows from these observations that $\mathcal{P}_r(v_2)=\mathcal{P}_r(e)$ 
if $e$ is the rightmost ingoing red edge at its end-vertex, 
that $\mathcal{P}_r(v_1)=\mathcal{P}_r(e)$ if $e$ is the leftmost outgoing 
red edge at its origin, and that no vertex $v$ satisfies 
$\mathcal{P}_r(v)=\mathcal{P}_r(e)$ otherwise.\hfill\phantom{1}
\end{proof}

\noindent\textbf{Definition.} Given an irreducible triangulation $T$
 endowed with its minimal transversal edge-partition, 
a \emph{ccw-internal} edge of $T$ is an inner edge $e$ of $T$ such 
that the counterclockwise-consecutive edge at each extremity of $e$ has the 
same color as $e$. Hence, on the minimal transversal pair of bipolar 
orientations, a ccw-internal red edge is an inner red edge of $T$ 
that is neither the leftmost outgoing red edge at its origin nor the 
rightmost ingoing red edge at its end-vertex.

\begin{lemma}
\label{prop:unused}
The number of abscissas not used by \textsc{Draw}($T$) is equal to 
the number of ccw-internal red edges of $T$. Similarly, the number of 
ordinates not used by \textsc{Draw}($T$) is equal to the number of 
ccw-internal blue edges of $T$.
\end{lemma}
\begin{proof}
Let $1\leq i\leq f_r$ be an abscissa-candidate and let $f$ be the unique 
inner face of $T_r$ such that $\mathrm{Abs}(f)=i$. Recall that the 
separating red path $\mathcal{P}_r(f)$ is equal to the separating red path 
of the bottom-right edge $e_f$ of $f$. Lemma~\ref{lemma:sepvertedge} ensures 
that $i$ is not the abscissa of any vertex of $T$ iff $e_f$ is a ccw-internal 
red edge of $T$. Hence, the number of abscissas not used by 
\textsc{Draw}($T$) is equal to the number of ccw-internal red edges 
of $T$ that are the bottom-right edge of an inner face of $T_r$. 
This quantity is also the number of ccw-internal red edges of $T$. Indeed, a 
ccw-internal red edge $e$ is the bottom-right edge of the inner face of $T_r$ 
on its left, as $e$ is not the leftmost outgoing red edge at its 
origin.\hfill\phantom{1}   
\end{proof}

\subsubsection{Reduction to the analysis of a parameter on ternary trees}
According to Lemma~\ref{prop:unused}, the number of coordinate-lines that
bear no vertex in \textsc{Draw}($T$) is equal to the 
number of ccw-internal edges of $T$. As we show next, 
the ccw-internal edges of $T$ correspond 
to particular edges of the ternary tree $A$ obtained by doing the 
opening of $T$. We define an \emph{internal edge} of a ternary 
tree $A$ as a closed edge $e$ of $A$ such that both edges 
following $e$ in clockwise order around each extremity of $e$ 
are closed edges.

\begin{lemma}
\label{lemma:internal_closure}
Let $T$ be an irreducible triangulation  and let $A$ be the 
ternary tree obtained by doing the opening of $T$, $T$ being endowed with 
its minimal transversal 
edge-partition and $A$ being endowed with the induced 
edge-bicoloration. Then each ccw-internal red (blue) edge of $T$ 
corresponds, via the opening of $T$, to an internal red (blue, respectively) 
edge of $A$. 
\end{lemma}  
\begin{proof}
Let $e=(v,v')$ be a ccw-internal red edge of $T$. By definition, the 
ccw-consecutive edge of $e$ at $v'$ is red. Hence the cw-consecutive 
edge of $e$ at $v$ is blue (otherwise there would be a unicolored face 
in $T$, contradicting Lemma~\ref{lemma:bic}). Similarly, the cw-consecutive 
edge of $e$ at $v'$ is blue. Hence, by definition of the opening mapping, $e$ 
is a closed edge of $A$. To prove that $e$ is an internal edge of $A$, we have to prove that 
the cw-consecutive edge of $A$ at each extremity of $e$ is a closed 
edge of $A$. Let $e_1,\ldots,e_k$ be the clockwise interval of blue edges of 
$T$ following $e$ in clockwise order around $v$ (the word interval refeering 
to the terminology of Condition C1); and let $v_1,\ldots,v_k$ be the 
corresponding sequence of neighbours of $v$. By definition of the opening 
mapping, the edge of $A$ following $e$ in clockwise order around $v$ 
is $e_k$. To prove that $e_k$ is a closed edge of $A$ it remains to show that 
the edge $e'\in T$ following $e_k$ in clockwise order around $v_k$ is red. 
If $k=1$, $e'$ is also the edge following $e$ in counterclockwise order 
around $v'$. Hence, the fact that $e$ is internal red ensures that $e'$ 
is red. If $k\geq 2$, then $e'=(v_k,v_{k-1})$. As $(v,v_k)$ and $(v,v_{k-1})$ 
are blue, Lemma~\ref{lemma:bic} ensures that $(v_k,v_{k-1})$ is red. 
Hence, $e_k$ is a closed edge of $A$. Similarly, the edge of $A$ following $e$ 
in clockwise order around $v'$ is a closed edge of $A$. Thus, $e$ is an internal edge of $A$.

Conversely, let $e=(v,v')$ be an internal red edge of $A$. 
By definition of the opening mapping, the edge $e_v\in T$
following $e$ in ccw order around $v$ has the same color as 
$e$ iff the half-edge of $e_v$ incident to $v$ has been created
during a local closure. The same holds with the edge $e_{v'}$ following
$e$ in ccw order around $v'$. We claim that this condition is satisfied by $e$ 
at $v'$ (and similarly at $v$). Indeed, the sides of closed edges incident to 
the outer face of the figure $F$ obtained 
as the partial closure of $A$ are such that the cw-consecutive 
edge at their right extremity (looking toward the outer face) 
is a stem. Hence, $e$ is not incident
to the outer face of $F$, so that there exists a local closure whose
effect is to close a triangular face incident to the right side of $e$
(traversed from $v$ to $v'$). The stem $s$ involved in this local closure
 can not be incident to $v$, as the edge of $A$ following $e$ in cw order 
around $v$ is not a stem. Hence $e=(v,v')$ is the second closed edge following 
the stem, so that the newly created half-edge (opposite to $s$) is incident
to $v'$. This concludes the proof.\hfill\phantom{1}
\end{proof}

Lemma~\ref{lemma:internal_closure} yields the following result that, as 
Fact~\ref{fact:abscissa}, follows from an easy adaptation of the proof of 
Theorem~\ref{cor:bije}.

\begin{fact}
\label{fact:disdeux}
The distribution of the parameter $\Delta_{\mathrm{abs}}(T)$, for $T$ taken 
uniformly at random in $\mathcal{T}_n'$, is equal to the distribution of the 
number of internal red edges in a rooted edge-bicolored ternary tree with $n$ 
 nodes taken uniformly at random.
\end{fact}
Let $X_n$ be the random variable denoting the number of internal red 
edges of a rooted edge-bicolored ternary tree with $n$  nodes taken 
uniformly at random. Fact~\ref{fact:disdeux} and the discussion in the 
overview ensure that the statement of 
Theorem~\ref{prop:reduction} about $W_c$ is proved by analyzing the 
distribution of $X_n$ and showing a concentration around $5n/54$.

\subsubsection{Analysis of the random variable $X_n$} 

\begin{figure}
\begin{center}
\includegraphics{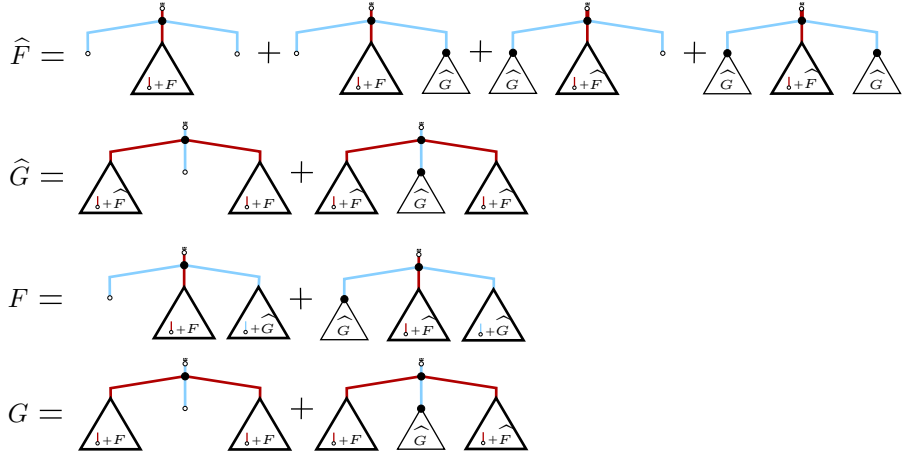}
\end{center}
\caption{Decomposition at the root node keeping track of the number of 
internal red edges.}
\label{fig:decomp_internal}
\end{figure}

We introduce the generating functions $F(z,u)$ and $G(z,u)$ counting 
respectively red-rooted ternary trees and blue-rooted ternary trees with 
respect to the number of  nodes and the number of internal red edges. 
As for the series $R(z,u)$ and $B(z,u)$ involved in the analysis of 
\textsc{Draw}, an equation system linking the series $F(z,u)$ 
and $G(z,u)$ can be obtained by decomposing a rooted ternary tree at its 
root node into three subtrees. 

To keep track of the parameter counting 
internal red edges, we introduce two auxiliary generating functions: 
$\widehat{F}(z,u)$ is the series counting red-rooted ternary trees with 
respect to the number of  nodes and the number of internal red edges, 
with  the difference that the root red stem is also considered as an internal 
red edge if its right child-edge is a closed edge; and $\widehat{G}(z,u)$ is 
the series counting blue-rooted ternary trees with respect to the number 
of  nodes and the number of internal red edges, with  the difference 
that the left child-edge $e$ of the root node is also considered as an 
internal red edge if $e$ is a closed edge and if the right child-edge of $e$ is a closed edge. 
The decomposition at the root node keeping track of the number of internal 
red edges is shown in Figure~\ref{fig:decomp_internal}. This directly 
translates to the following system:
\begin{equation}
\left\{
\begin{array}{rcl}
\widehat{F}(z,u) &=& z(1+F)+zu(1+F)\widehat{G}+z(1+\widehat{F})\widehat{G}+zu\widehat{G}(1+\widehat{F})\widehat{G}\\
\widehat{G}(z,u) &=& z(1+\widehat{F})(1+F)+z(1+\widehat{F})\widehat{G}(1+\widehat{F})\\
F(z,u) &=& z(1+F)(1+\widehat{G})+z\widehat{G}(1+\widehat{F})(1+\widehat{G})\\
G(z,u) &=& z(1+F)^2+z(1+F)\widehat{G}(1+\widehat{F})
\end{array}
\right.
\label{eq:FG}
\end{equation}

Let $H(z,u)$ be the series counting rooted edge-bicolored ternary trees 
by the numbers of nodes and internal red edges; 
clearly $H(z,u)=F(z,u)+G(z,u)$. 
An algebraic equation $P(H(z,u),z,u)=0$ is easily
computed from~(\ref{eq:FG}). Then, the algebraic quasi-power theorem 
can be applied to the algebraic generating function $H(z,u)$. All conditions 
are easily checked and the algebraic series $\rho(u)$ which we obtain 
satisfies $\mu=-\rho'(1)/\rho(1)=5/54$. This yields the statement of 
Theorem~\ref{prop:reduction} on the distribution of $W_c$. 
The result is the same for $H_c$, because the distribution of $H_c$ 
is equal to the distribution of $W_c$, by stability of $\mathcal{T}_n'$ 
and of \textsc{CompactDraw} under the $\pi/2$ clockwise rotation. 
This concludes the proof of 
Theorem~\ref{prop:reduction}.
%%???%%\vspace{-0.4em}

\section{Conclusion}
Several families of planar maps are characterised by the existence
of a specific combinatorial structure for each map of the family
(bipolar orientations for 2-connected maps, Schnyder woods for 
triangulations). This article provides a detailed study of 
so-called \emph{transversal structures}, which are specific 
to the family of triangulations of the 4-gon with no filled 3-cycle ---called
irreducible triangulations. The two main results we obtain related
to transversal structures are a bijection between irreducible triangulations
and ternary trees, as well as a new straight-line drawing algorithm
for irreducible triangulations. The bijection allows us to count 
irreducible triangulations and the closely related family of 4-connected
triangulations. 

Surprisingly, the bijection ---combined with some modern 
tools of analytic combinatorics--- also makes it possible to analyse
the grid size of the straight-line drawing; for a random irreducible
triangulation with $n$ vertices, the grid size is with high 
probability $11n/27\times 11n/27$ up to an additive error $\cO(\sqrt{n})$.
In comparison, the best previously known algorithm~\cite{Miura} only
guarantees a grid size $(\lceil n/2\rceil -1)\times \lfloor n/2\rfloor$.
The principle of our straight-line drawing algorithm has also been 
recently applied to the family of quadrangulations~\cite{Fu06,Fu-these}; for
a random quadrangulation with $n$ vertices, the grid
size is with high probability $13n/27\times 13n/27$ up to an additive
error $\cO(\sqrt{n})$, which improves on the grid size 
$(\lceil n/2\rceil -1)\times \lfloor n/2\rfloor$ guaranteed by~\cite{Bi}.

\vspace{0.2cm}

%%%%%%%%%%%%%%%%%%%%%%%%%%%%%%%%%%%%%%%%%%%%%%%%%%%%%%%%%%%%%%%%%%%%%%%%%%%
\noindent{\bf Acknowledgments.} I would like to thank my advisor Gilles
Schaeffer, who has greatly
helped me to produce this work
through numerous discussions, steady encouragment and useful suggestions. 
I also thank Nicolas Bonichon, Luca Castelli Aleardi, and Philippe Flajolet 
for fruitful discussions, Roberto Tamassia for having pointed
reference~\cite{To} to me, and Thomas Pillot for 
very efficient implementations of all algorithms presented in this article. 
Finally, I am very grateful to the two anonymous referees
for their detailed corrections and insightful remarks.

\bibliographystyle{plain}
\bibliography{DiscreteMathFusy}  

\newpage
\appendix

\end{document}